\documentclass[10pt, twoside]{article}
\usepackage{amsmath,amsthm,amssymb}
\usepackage{times}
\usepackage{enumerate}
\usepackage{amssymb}
\usepackage{graphicx}
\usepackage[noadjust]{cite}

\usepackage[dvipsnames]{color}

\def\titlerunning#1{\gdef\titrun{#1}}
\makeatletter
\def\author#1{\gdef\autrun{\def\and{\unskip, }#1}\gdef\@author{#1}}

\makeatother

\def\keywords#1{\par\medskip
\noindent\textbf{Keywords.} #1}
\def\subjclass#1{\par\smallskip
\noindent\textbf{MSC (2010):} #1}

\newtheorem{thm}{Theorem}[section]

\newtheorem{lem}[thm]{Lemma}

\theoremstyle{definition}

\newtheorem{rem}[thm]{Remark}
\newtheorem{exa}[thm]{Example}

\numberwithin{equation}{section}

\newtheorem*{notations}{Notations}

\DeclareMathOperator*{\esssup}{ess\,sup}
\DeclareMathOperator*{\essinf}{ess\,inf}

\makeatletter
\let\@fnsymbol\@alph
\makeatother

\begin{document}

\baselineskip=17pt

\titlerunning{Radial quasilinear problems}

\title{Radial quasilinear elliptic problems with singular or vanishing potentials}

\author{Marino Badiale\thanks{Dipartimento di Matematica ``Giuseppe Peano'', Universit\`{a} degli Studi di
Torino, Via Carlo Alberto 10, 10123 Torino, Italy. 
e-mails: \texttt{marino.badiale@unito.it}, \texttt{michela.guida@unito.it}}
\textsuperscript{,}\thanks{Partially supported by the PRIN2012 grant ``Aspetti variazionali e
perturbativi nei problemi di.renziali nonlineari''.}
\ -\ Michela Guida\textsuperscript{a}
\ -\ Sergio Rolando\thanks{Dipartimento di Matematica e Applicazioni, Universit\`{a} di Milano-Bicocca,
Via Roberto Cozzi 53, 20125 Milano, Italy. e-mail: \texttt{sergio.rolando@unito.it}}%
}

\date{
}
\maketitle

\begin{abstract}
In this paper we continue the work that we began in \cite{AVK_I}. 
Given $1<p<N$, two measurable functions $V\left(r \right)\geq 0$ and $K\left(r\right)> 0$, 
and a continuous function $A(r) >0$ ($r>0$), we consider the quasilinear elliptic equation
\[
-\mathrm{div}\left(A(|x| )|\nabla u|^{p-2} \nabla u\right) +V\left( \left| x\right| \right) |u|^{p-2}u= K(|x|) f(u) \quad \text{in }\mathbb{R}^{N},
\]
where all the potentials $A,V,K$ may be singular or vanishing, at the origin or at infinity. We find existence of nonnegative solutions by the application of variational methods, for which we need to study the compactness of the embedding of a suitable
function space $X$ into the sum of Lebesgue spaces $L_{K}^{q_{1}}+L_{K}^{q_{2}}$. 
The nonlinearity has a double-power super $p$-linear behavior, as $f(t)= 
\min \left\{ t^{q_1 -1}, t^{q_2 -1}     \right\}$ with $q_1,q_2>p$ (recovering the power case if $q_1=q_2$). With respect to \cite{AVK_I}, in the present paper we assume some more hypotheses on $V$, and we are able to enlarge the set of values $q_1 , q_2$ for which we get existence results.

\end{abstract}

\section{Introduction}

In this paper, we pursue the work we made in \cite{BGrR, BGR_I,BGR_II, BGR_p, AVK_I, BZ, GR-nls} and, in particular, we complete the results we proved in \cite{AVK_I}. 
In all these works, we have obtained embedding and compactness results for weighted Sobolev spaces, which, by variational methods, yielded existence and multiplicity results for nonlinear elliptic equations in $\mathbb{R}^N$. 
In \cite{AVK_I} and in the present paper, we consider quasilinear elliptic equations in presence of a radial potential on the derivatives, namely, equations of the following form:
\begin{equation}\label{EQ}
-\mathrm{div}\left(A(|x|) |\nabla u|^{p-2} \nabla u\right) +V\left( \left| x\right| \right) |u|^{p-2}u= K(|x|) f(u) \quad \text{in }\mathbb{R}^{N},
\end{equation}
where $1<p<N$, $f:\mathbb{R}\rightarrow \mathbb{R}$ is a continuous nonlinearity satisfying $f\left( 0\right) =0$, and $V\geq 0$ and $A,K>0$ are given potentials, which may be unbounded or vanishing at the origin or at infinity. 
This equation has been recently studied by Su and Wang in \cite{Su-Wang} (see also the references therein), where the authors considered potentials $A(r)$, $V(r)$ and $K(r)$ essentially behaving as powers of $r$ as $r\to 0$ and $r\to +\infty$ (see the paper introduction for more precise assumptions), and investigated the effects on the solutions structure from the interplay between the quasilinearity $p$, the nonlinear growth of $f(u)$ and the growth or decay rates of the potentials. The aim of the present paper is to extend the results of \cite{Su-Wang} (see below for further detalis).

We study Equation \eqref{EQ} by variational methods, so we introduce a suitable functional space $X$ (see section 2) and we say that $u\in X$ is a \textit{weak solution}\emph{\ }to (\ref{EQ}) if 
\begin{equation*}
\int_{\mathbb{R}^{N}}A(|x|)|\nabla u|^{p-2}\nabla u\cdot \nabla h\,dx+\int_{\mathbb{R}^{N}}V\left(
\left| x\right| \right) |u|^{p-2}uh\,dx=\int_{\mathbb{R}^{N}}K\left( \left| x\right|
\right) f\left( u\right) h\,dx
\end{equation*}
\noindent for all $h \in X$. These solutions are (at least formally) critical points of the Euler
functional 
\begin{equation}
\label{I:=}
I\left( u\right) :=\frac{1}{p}\left\| u\right\| ^{p}-\int_{\mathbb{R}
^{N}}K\left( \left| x\right| \right) F\left( u\right) dx
\end{equation}
where $||\cdot ||$ is the norm in $X$, so the problem of existence is quite easily solved by the use of standard embedding results for Sobolev spaces if $A$ is bounded and bounded away from zero, $V$ does not vanish at infinity and $K$ is bounded. These results are not available anymore for potentials which may be singular or vanishing, at zero or at infinity, and new embedding theorems are needed. 
Much work has been done in recent years in this direction, and we refer the reader for example to
the following papers and the references therein: 
\cite{BGR_I,BGR_II,GR-nls,Su-Wang-Will-2} for the case $A=1,p=2$ of the usual laplacian operator, 
\cite{Anoop,Su12,Cai-Su-Sun,SuTian12,Su-Wang-Will-p,Yang-Zhang,Zhang13,BPR, BGR_p} for the $p$-laplacian case $A=1,1<p<N$  
(see also \cite{BGrR} for equations involving the bi-laplacian operators), 
\cite{BZ} for the case $A\neq 1, p=2$ of the weighted laplacian, and the already mentioned paper \cite{Su-Wang}, which inspired the present work, for the general case $A\neq 1,1<p<N$ with a potential on the derivatives. 

The main novelty of our approach with respect to \cite{Su-Wang} (here and in \cite{AVK_I}) is two-fold. First, we look for embeddings of $X$ not into a single Lebesgue space $L_{K}^{q}$ 
but into a sum of Lebesgue spaces $L_{K}^{q_1}+L_{K}^{q_2}$. This allows us to study separately the behaviour of the potentials $V$ and $K$ at $0$ and $\infty$, and thus to assume independent sets of hypotheses about these behaviours, which therefore do not need to be linked to one another to get existence. Second, we assume hypotheses not on $V$ and $K$ separately but on their ratio, so admitting asymptotic behaviors of general kind for the two potentials (not only power-like).
In particular, with respect to \cite{AVK_I}, here we assume some more specific hypotheses on the asymptotic behavior of the potential $V$, and this allows us to get larger intervals of exponents $q_1 ,q_2$ for which the embedding 
$X\hookrightarrow L_{K}^{q_1}+L_{K}^{q_2}$ is compact, and then equation (\ref{EQ}) has solutions. The main results of this paper are Theorems \ref{THM2} and \ref{THM3} (about compact embeddings), and Theorem \ref{THM:ex2} (about existence of solutions).

As a conclusion, putting together the results of \cite{AVK_I} and the present paper, we extend the results of \cite{Su-Wang} by getting existence results for large new classes of nonlinearities and radial potentials. This is exemplified, both here and in \cite{AVK_I}, in a specific section of examples, where we show existence for some particular cases for which nothing was known in the previous literature (see Section \ref{SEC:EX} below and Section 7 in \cite{AVK_I}).

The paper is organized as follows. In Section \ref{SEC:MAIN} we introduce our hypotheses on $A,V,K$ 
and the function spaces $X$ and $L_{K}^{q_{1}}+L_{K}^{q_{2}}$ in which we will work. 
In section \ref{COMP} we give a general result concerning the embedding properties of $X$ into $L_{K}^{q_{1}}+L_{K}^{q_{2}}$ (Theorem \ref{THM(cpt)}),
and then we state some explicit conditions ensuring that such an embedding is compact (Theorems \ref{THM2} and \ref{THM3}). 
These compactness theorems are proved in section \ref{SEC:2}. 
In Section \ref{SEC: ex}, we apply our previous results to get existence of nonnegative solutions for (\ref{EQ}). The main existence theorem of this paper is Theorem \ref{THM:ex2}. In section \ref{SEC:EX}, we give some examples of specific cases in which our results apply, in order to explain their novelties with respect to the literature. Finally, the Appendix is devoted to some detailed computations, displaced there from Section \ref{SEC:2} for sake of clarity.

\begin{notations}

We end this introductory section by collecting
some notations used in the paper.

\noindent $\bullet $ For every $R>0$, we set $B_{R}:=\left\{ x\in \mathbb{R}
^{N}:\left| x\right| <r\right\} $.

\noindent $\bullet $ For any subset $A\subseteq \mathbb{R}^{N}$, we denote $
A^{c}:=\mathbb{R}^{N}\setminus A$. If $A$ is Lebesgue measurable, $\left|
A\right| $ stands for its measure.



\noindent $\bullet $ $C_{\mathrm{c}}^{\infty }(\Omega )$ is the space of the
infinitely differentiable real functions with compact support in the open
set $\Omega \subseteq \mathbb{R}^{N}$. If $\Omega$ has radial symmetry, $C_{\mathrm{c}, r}^{\infty }( \Omega )$ is the subspace 
of $C_{\mathrm{c}}^{\infty }(\Omega )$ made of radial functions.

\noindent $\bullet $ If $1\leq p\leq \infty $ then $L^{p}(A)$ and $L_{
\mathrm{loc}}^{p}(A)$ are the usual real Lebesgue spaces (for any measurable
set $A\subseteq \mathbb{R}^{N}$). If $\rho :A\rightarrow \left( 0,+\infty
\right) $ is a measurable function, then $L^{p}(A,\rho \left( z\right) dz)$
is the real Lebesgue space with respect to the measure $\rho \left( z\right)
dz$ ($dz$ stands for the Lebesgue measure on $\mathbb{R}^{N}$).

\noindent $\bullet $ $p^{\prime }:=p/(p-1)$ is the H\"{o}lder-conjugate
exponent of $p.$

\noindent  $\bullet $ $\mathbb{R}_{+} = ( 0, +\infty )$.

\end{notations}

\section{Hypotheses and preliminary results \label{SEC:MAIN} }

Throughout the paper we assume $N\geq 3$, $1<p<N$ and the following set of hypotheses on $A, V,K$:

\begin{itemize}

\item[$\left( \mathbf{A}\right) $]  $A : \mathbb{R}_{+} \rightarrow \mathbb{R}_{+} $ is a continuous function 
and there exist $a_{0},a_{\infty} \in(p-N,p]$ and $c_{0},c_{\infty}>0$ such that:
$$
c_{0} \leq \liminf_{r\rightarrow 0^+}\frac{A(r)}{r^{a_{0}}} \leq \limsup_{r\rightarrow 0^+}\frac{A(r)}{r^{a_{0}}} < + \infty, 
$$
$$
c_{\infty}\leq \liminf_{r\rightarrow +\infty}\frac{A(r)}{r^{a_{\infty}}} \leq \limsup_{r\rightarrow +\infty}\frac{A(r)}{r^{a_{\infty}}} <+\infty 
$$

\item[$\left( \mathbf{V}\right) $]  $V:\mathbb{R}_{+} \rightarrow
\left[ 0,+\infty \right) $ is a measurable function such that $V\in
L_{\mathrm{loc}}^{1}\left( \mathbb{R}_{+}  \right) $

\item[$\left( \mathbf{K}\right) $]  $K:\mathbb{R}_{+} \rightarrow \mathbb{R}_{+} $ is a measurable function such that $K\in L_{
\mathrm{loc}}^{s}\left( \mathbb{R}_{+}  \right) $ for some $s>1$.

\end{itemize}







Under assumption $\left( \mathbf{A}\right) $, we define
\[
S_A := \left\{  u \in C_{\mathrm{c}, r}^{\infty }( \mathbb{R}^{N} )  \, \Big| \, \int_{\mathbb{R}^{N}}A(|x|) |\nabla u |^p dx <+\infty  \right\}
\quad\text{and}\quad
||u||_A := \left( \int_{\mathbb{R}^{N}}A(|x|) |\nabla u |^p dx \right)^{1/p},
\]
and denote by $D_A$ the completion of $S_A$ with respect to the norm $|| \cdot ||_A$.
The main properties of $D_A$ are summarized in Lemma \ref{lemdens2} below (see \cite{AVK_I} for a proof), where 
the exponents $p_{0}, p_{\infty}\geq p$ are defined by
\begin{equation*}
p_{0}:=\frac{pN}{N+a_{0}-p}\quad \text{and}\quad p_{\infty}:=\frac{pN}{N+a_{\infty}-p},
\end{equation*}
Notice that $p_{0}, p_{\infty}\geq p$ thanks to the hypothesis $a_{0},a_{\infty} \in(p-N,p]$. 
The main function space we will work in is then defined by setting
\[
X := D_A \cap  L^{p}(\mathbb{R}^{N}, V (|x|) dx ) \quad\text{and}\quad ||u|| := \left( ||u||_{A}^{p} + ||u||_{L^{p}(\mathbb{R}^{N}, V (|x|) dx )}^{p} 
\right)^{1/p}. 
\]
Note that $X$ is nonempty by assumption $\left( \mathbf{A}\right) $ and it is a Banach space with respect to the norm $||\cdot||$.

\begin{lem}
	\label{lemdens2}
	The following properties hold.
	
	\begin{itemize}

\item[{$(i)$}] If $u \in D_A$, then $u$ has weak derivatives $D_i u$ in the open set $\Omega=\mathbb{R}^{N} \backslash \{0 \}$ ($i=1,..., N$) and $D_i u \in L^p_{loc}(\Omega )$.

\item[{$(ii)$}] If $u \in D_A$ then 
$||u||_A =\left( \int_{\mathbb{R}^{N}} A(|x|) |\nabla u|^p \, dx \right)^{1/p}$ is finite and defines a norm on $D_A$. With this norm, $D_A$ is a Banach space.

\item[{$(iii)$}] For any $R_0 >0$, there exists a constant $C=C(N, R_0 , a_0 , a_{\infty})>0$ such that for all $u \in D_A$ one has
$$ |u(x) | \leq C |x|^{-\frac{N+ a_0 -p}{p}} \, ||u||_A \quad \text{almost everywhere in } B_{R_0}, $$
$$ |u(x) | \leq C |x|^{-\frac{N+ a_{\infty} -p}{p}} \, ||u||_A \quad \text{almost everywhere in } B^{c}_{R_0}. $$

\item[{$(iv)$}] For any $0<r<R$, 
the embeddings $D_A \hookrightarrow L^{p_{0}}(B_{R})$ and $D_A \hookrightarrow L^{p_{\infty}}(B^{c}_{R})$
are continuous, and the embedding $D_{A} \hookrightarrow L^{p}(B_R \backslash {\overline B_r })$
is continuous and compact.

\end{itemize}

\end{lem}

We end this section by recalling some definitions and results from \cite{BPR} concerning
the other function space which will be relevant in the following, i.e., the sum space 
\[
L_{K}^{q_{1}}+L_{K}^{q_{2}}:=L_{K}^{p_{1}}\left( \mathbb{R}^{N}\right)
+L_{K}^{q_{2}}\left( \mathbb{R}^{N}\right) :=\left\{ u_{1}+u_{2}:u_{1}\in
L_{K}^{q_{1}}\left( \mathbb{R}^{N}\right) ,\,u_{2}\in L_{K}^{q_{2}}\left( \mathbb{R
}^{N}\right) \right\} , 
\]
where we assume $1<q_{1}\leq q_{2}<\infty $. Such a space can be
characterized as the set of the measurable mappings $u:\mathbb{R}
^{N}\rightarrow \mathbb{R}$ for which there exists a measurable set $E\subseteq 
\mathbb{R}^{N}$ such that $u\in L_{K}^{q_{1}}\left( E\right) \cap
L_{K}^{q_{2}}\left( E^{c}\right) $ (of course $L_{K}^{q_{1}}\left( E\right)
:=L^{q_{1}}(E,K\left( \left| x\right| \right) dx)$, and so for $
L_{K}^{q_{2}}\left( E^{c}\right) $). It is a Banach space with respect to
the norm 
\[
\left\| u\right\| _{L_{K}^{q_{1}}+L_{K}^{q_{2}}}:=\inf_{u_{1}+u_{2}=u}\max
\left\{ \left\| u_{1}\right\| _{L_{K}^{q_{1}}},\left\| u_{2}\right\|
_{L_{K}^{q_{2}}}\right\} 
\]
and the continuous embedding $L_{K}^{q}\hookrightarrow
L_{K}^{q_{1}}+L_{K}^{q_{2}}$ holds for all $q\in \left[ q_{1},q_{2}\right] $.



\section{Compactness results}\label{COMP}

In this section we state the main compactness results of this paper. 
Recall that we let $N \geq 3$, $1<p<N$ and $A$, $V$, $K$ be as in 
$\left( \mathbf{A}\right) $, $\left( \mathbf{V}\right) $, $\left( \mathbf{K}\right) $.

We define the following functions of $R>0$ and $q>1$: 
\begin{eqnarray}
\mathcal{S}_{0}\left( q,R\right)&:=&
\sup_
{u\in X,\,
\left\| u\right\| =1  }
\int_{B_{R}}K\left( \left| x\right| \right)
\left| u\right| ^{q}dx,  \label{S_o :=}
\\
\mathcal{S}_{\infty }\left( q,R\right)&:=&
\sup_
{u\in X,\,
\left\| u\right\| =1  } 
\int_{\mathbb{R}%
^{N}\setminus B_{R}}K\left( \left| x\right| \right) \left| u\right| ^{q}dx.
\label{S_i :=}
\end{eqnarray}
Clearly $\mathcal{S}_{0}\left( q,\cdot \right) $ is nondecreasing, $\mathcal{
S}_{\infty }\left( q,\cdot \right) $ is nonincreasing and both of them can
be infinite at some $R$.

Our first result (Theorem \ref{THM(cpt)}) concerns embedding properties of $X$
into $L_{K}^{q_{1}}+L_{K}^{q_{2}}$ and relies on assumptions which are quite
general, but not so easy to check.
More handy conditions ensuring these general assumptions will be provided then by Theorems \ref{THM2} and \ref{THM3}.

\begin{thm}
\label{THM(cpt)} 
Let $q_{1},q_{2}>1$.

\begin{itemize}
\item[(i)]  If 
\begin{equation}
\mathcal{S}_{0}\left( q_{1},R_{1}\right) <\infty \quad \text{and}\quad 
\mathcal{S}_{\infty }\left( q_{2},R_{2}\right) <\infty \quad \text{for some }
R_{1},R_{2}>0,  
\tag*{$\left( {\cal S}_{q_{1},q_{2}}^{\prime }\right) $}
\end{equation}
then $X$ is continuously embedded into $L_{K}^{q_{1}}(\mathbb{R}^{N})+L_{K}^{q_{2}}(\mathbb{R}^{N})$.

\item[(ii)]  If 
\begin{equation}
\lim_{R\rightarrow 0^{+}}\mathcal{S}_{0}\left( q_{1},R\right)
=\lim_{R\rightarrow +\infty }\mathcal{S}_{\infty }\left( q_{2},R\right) =0, 
\tag*{$\left({\cal S}_{q_{1},q_{2}}^{\prime \prime }\right) $}
\end{equation}
then $X$ is compactly embedded into $L_{K}^{q_{1}}(\mathbb{R}^{N})+L_{K}^{q_{2}}(\mathbb{R}^{N})$.

\end{itemize}
\end{thm}

\noindent For the proof of Theorem \ref{THM(cpt)} see \cite{AVK_I}. Observe that, of course, $(\mathcal{S}_{q_{1},q_{2}}^{\prime \prime })$
implies $(\mathcal{S}_{q_{1},q_{2}}^{\prime })$. Moreover, these assumptions
can hold with $q_{1}=q_{2}=q$ and therefore Theorem \ref{THM(cpt)} also
concerns embedding properties of $X$ into $L_{K}^{q}$, $1<q<\infty $.

Explicit conditions implying $(\mathcal{S}_{q_{1},q_{2}}^{\prime })$ for some $q_{1},q_{2}$ have been given in Theorems 3.2 and 3.3 of \cite{AVK_I}. 
Here we pursue this work and in Theorems \ref{THM2} and \ref{THM3} below we improve the results of \cite{AVK_I} by exploiting further informations on the growth of $V$. 
More precisely, in Theorem \ref{THM3} we will find a range of exponents $
q_{1} $ such that $\lim_{R\rightarrow 0^{+}}\mathcal{R}_{0}\left(q_{1},R\right)$ $=0$, 
while in Theorem \ref{THM2} we will do the same for exponents $q_{2}$ such that
$\lim_{R\rightarrow +\infty}\mathcal{R}_{\infty }\left( q_{2},R\right) =0$.
Condition $(\mathcal{S}_{q_{1},q_{2}}^{\prime \prime })$ 
then follows by joining Theorem \ref{THM3} with Theorem \ref{THM2}.

To state Theorems \ref{THM2} and \ref{THM3}, we need the following functions. For any $a>p-N$, $\alpha \in \mathbb{R}$, $\beta \leq 1$ and $\gamma \in \mathbb{R}$, define 
\begin{equation}
q_{*}\left( \alpha ,\beta ,\gamma \right) :=p\, \frac{\alpha -\gamma \beta +N}{%
N-\gamma }\quad \text{and}\quad q_{**}\left(a, \alpha ,\beta ,\gamma \right)
:=
p\, \frac{p\alpha +\left( 1-p\beta \right) \gamma +p\left( N-1\right)+a }{%
p\left( N-1\right) -\gamma (p-1) +a}.  \label{q** :=}
\end{equation}
Of course $q_{*}$ and $q_{**}$ are undefined if $\gamma =N$ and $\gamma
= \frac{p (N-1) +a}{p-1} $, respectively

\begin{thm}
\label{THM2}
Assume that there exists $R_{2}>0$ such that
\begin{equation}
\esssup_{r>R_{2}}\frac{K\left( r\right) }{r^{\alpha _{\infty
}}V\left( r\right) ^{\beta _{\infty }}}<+\infty \quad \text{for some }0\leq
\beta _{\infty }\leq 1\text{~and }\alpha _{\infty }\in \mathbb{R}
\label{hp all'inf}
\end{equation}
and 
\begin{equation}
\essinf_{r>R_{2}}r^{\gamma _{\infty }}V\left( r\right) >0\quad 
\text{for some }\gamma _{\infty } \leq p - a_{\infty}.  \label{stima all'inf}
\end{equation}
Then $\displaystyle \lim_{R\rightarrow +\infty }\mathcal{S}_{\infty }\left(
q_{2},R\right) =0$ for every $q_{2}\in \mathbb{R}$ such that 
\begin{equation}
q_{2}>\max \left\{ 1,p\beta _{\infty },q_{*},q_{**}\right\} ,  \label{th3}
\end{equation}
where $q_{*}=q_{*}\left( \alpha _{\infty },\beta _{\infty },\gamma _{\infty
}\right) $ and $q_{**}=q_{**}\left(  a_{\infty }, \alpha _{\infty },\beta _{\infty
},\gamma _{\infty }\right) .$
\end{thm}

\begin{rem}
The proof of Theorem \ref{THM2} does not
require $\beta _{\infty }\geq 0$, but this condition is not a restriction of
generality in stating the theorem. Indeed, under assumption (\ref{stima
all'inf}), if (\ref{hp all'inf}) holds with $\beta _{\infty }<0$, then it
also holds with $\alpha _{\infty }$ and $\beta _{\infty }$ replaced by $%
\alpha _{\infty }-\beta _{\infty }\gamma _{\infty }$ and $0$ respectively,
and this does not change the thesis (\ref{th3}), because $q_{*}\left( \alpha
_{\infty }-\beta _{\infty }\gamma _{\infty },0,\gamma _{\infty }\right)
=q_{*}\left( \alpha _{\infty },\beta _{\infty },\gamma _{\infty }\right) $
and $q_{**}\left( a_{\infty }, \alpha _{\infty }-\beta _{\infty }\gamma _{\infty
},0,\gamma _{\infty }\right) =q_{**}\left(  a_{\infty }, \alpha _{\infty },\beta _{\infty
},\gamma _{\infty }\right) $.
\end{rem}

In order to state our next result, we introduce, by the following
definitions, an open region $\mathcal{A}_{a, \beta ,\gamma }$ of the $\alpha q$
-plane, depending on $a>p-N$, $\beta\in[0,1]$ and $\gamma \geq p -a$. 
We set 
\begin{equation}
\begin{array}{ll}
\mathcal{A}_{a, \beta ,\gamma }:=\left\{ \left( \alpha ,q\right) :\max \left\{
1,p\beta \right\} <q<\min \left\{ q_{*},q_{**}\right\} \right\} \quad
\smallskip & \text{if }p-a\leq\gamma <N, \\ 
\mathcal{A}_{a,\beta ,\gamma }:=\left\{ \left( \alpha ,q\right) :\max \left\{
1,p\beta \right\} <q<q_{**},\,\alpha >-\left( 1-\beta \right) N\right\}
\quad \smallskip & \text{if }\gamma =N, \\ 
\mathcal{A}_{a,\beta ,\gamma }:=\left\{ \left( \alpha ,q\right) :\max \left\{
1,p\beta ,q_{*}\right\} <q<q_{**}\right\} \smallskip & \text{if }N<\gamma
<\frac{p(N-1)+a}{p-1}, \\ 
\mathcal{A}_{a, \beta ,\gamma }:=\left\{ \left( \alpha ,q\right) :\max \left\{
1,p\beta ,q_{*}\right\} <q,\,\alpha >-\left( 1-\beta \right) \gamma \right\}
\smallskip & \text{if }\gamma =\frac{p(N-1)+a}{p-1}, \\ 
\mathcal{A}_{a,\beta ,\gamma }:=\left\{ \left( \alpha ,q\right) :\max \left\{
1,p\beta ,q_{*},q_{**}\right\} <q\right\} & \text{if }\gamma >\frac{p(N-1)+a}{p-1}.
\end{array}
\label{A:=}
\end{equation}
Notice that $\frac{p(N-1)+a}{p-1}>N$ because $N-p+a>0$.

\begin{thm}
\label{THM3}
Assume that there exists $R_{1}>0$ such that 
\begin{equation}
\esssup_{r\in \left( 0,R_{1}\right) }\frac{K\left( r\right) }{%
r^{\alpha _{0}}V\left( r\right) ^{\beta _{0}}}<+\infty \quad \text{for some }%
0\leq \beta _{0}\leq 1\text{~and }\alpha _{0}\in \mathbb{R}  \label{hp in 0}
\end{equation}
and 
\begin{equation}
\essinf_{r\in \left( 0,R_{1}\right) }r^{\gamma _{0}}V\left(
r\right) >0\quad \text{for some }\gamma _{0}\geq p -a_{0}.
  \label{stima in 0}
\end{equation}
Then $\displaystyle \lim_{R\rightarrow 0^{+}}\mathcal{S}_{0}\left(
q_{1},R\right) =0$ for every $q_{1}\in \mathbb{R}$ such that 
\begin{equation}
\left( \alpha _{0},q_{1}\right) \in \mathcal{A}_{a_{0},\beta _{0},\gamma _{0}}.
\end{equation}
\end{thm}

\begin{rem}
\label{RMK: suff12}

 \label{RMK: suff12-V^0} In both Theorems \ref{THM2} and \ref{THM3}, we mean $V\left( r\right) ^{0}=1$ for every $r$
(even if $V\left( r\right) =0$). 

\end{rem}


\section{Proof of Theorems \ref{THM2} and \ref{THM3} \label
{SEC:2}}
This section is devoted to the proof of Theorems \ref{THM2} and \ref{THM3}, which will be achieved through several lemmas. Recall that we assume $N\geq 3$, $1<p<N$ and let $A$, $V$ and $K$ be as in $\left( \mathbf{A}\right)$, $\left( \mathbf{V}\right)$ and $\left( \mathbf{K}\right)$.

The first two lemmas are given in \cite{AVK_I} and we recall them here for completeness. 

\begin{lem}
\label{Lem(Omega)}Let $R_1>0$ and assume that 
\[
\Lambda_1 :=\esssup_{x\in B_{R_1} }\frac{K\left( \left| x\right|
\right) }{\left| x\right| ^{\alpha }V\left( \left| x\right| \right) ^{\beta }
}<+\infty \quad \text{for some }0\leq \beta \leq 1\text{~and }\alpha \in 
\mathbb{R}.
\]
Let $u\in X$ and assume that there exist $\nu \in \mathbb{R}$ and $m>0$
such that 
\[
\left| u\left( x\right) \right| \leq \frac{m}{\left| x\right| ^{\nu }}\quad 
\text{almost everywhere on } B_{R_1}.
\]
Then there exists a constant $C_1=C_1(N, p, a_0, R_1, \beta )>0$ such that $\forall R \in (0,R_{1})$ and $\forall q>\max \left\{ 1,p\beta \right\} $ one has 

$\displaystyle\int_{B_R }K\left( \left| x\right| \right) \left| u\right| ^q dx$
\[
\leq \left\{ 
\begin{array}{ll}
\Lambda_1 m^{q-1} C_1 \left( \int_{B_R}\left| x\right| ^{\frac{%
\alpha -\nu \left( q-1\right) }{N(p-1)+p\left( 1-p\beta +a_0 \beta \right) -a_0}pN}dx\right)^{
\frac{N(p-1)+p\left( 1-p\beta +a_0 \beta \right) -a_0}{pN}}\left\| u\right\| \quad \medskip  & 
\text{if }0\leq \beta \leq \frac{1}{p} \\ 
\Lambda_1 m^{q-p\beta }\left( \int_{B_R }\left| x\right| ^{\frac{\alpha
-\nu \left( q-p\beta \right) }{1-\beta }}dx\right) ^{1-\beta }\left\|
u\right\| ^{p\beta} \medskip  & \text{if }\frac{1}{p}
<\beta <1 \\ 
\Lambda_1 m^{q-p}\left( \int_{B_R }\left| x\right| ^{\frac{p}{p-1}(\alpha -\nu \left(
q-p\right)) }V\left( \left| x\right| \right) \left| u\right| ^{p}dx\right) ^{
\frac{p-1}{p}}\left\| u\right\|  & \text{if }\beta =1.
\end{array}
\right. 
\]
\end{lem}

\bigskip

\begin{lem}
\label{Lem(Omega2)} Let $R_2>0$ and assume that 
\[
\Lambda_2 :=\esssup_{x\in B_{R_2}^c }\frac{K\left( \left| x\right|
\right) }{\left| x\right| ^{\alpha }V\left( \left| x\right| \right) ^{\beta }%
}<+\infty \quad \text{for some }0\leq \beta \leq 1\text{~and }\alpha \in 
\mathbb{R}.
\]
Let $u\in X$ and assume that there exist $\nu \in \mathbb{R}$ and $m>0$
such that 
\[
\left| u\left( x\right) \right| \leq \frac{m}{\left| x\right| ^{\nu }}\quad 
\text{almost everywhere on }B_{R_2}^c .
\]
Then there exists a constant $C_2=C_2(N, p, a_{\infty}, R_2, \beta )>0$ such that $\forall R > R_{2}$ and $\forall q>\max \left\{ 1,p\beta \right\} $ one has 

$\displaystyle\int_{B_R^c }K\left( \left| x\right| \right) \left| u\right| ^{q}dx$
\[
\leq \left\{ 
\begin{array}{ll}
\Lambda_2 m^{q-1} C_2 \left( \int_{B_{R}^c}\left| x\right| ^{\frac{%
\alpha -\nu \left( q-1\right) }{N(p-1)+p\left( 1-p\beta +a_{\infty} \beta \right) -a_{\infty}}pN}dx\right)^{
\frac{N(p-1)+p\left( 1-p\beta +a_{\infty} \beta \right)-a_{\infty} }{pN}}\left\| u\right\| \quad \medskip  & 
\text{if }0\leq \beta \leq \frac{1}{p} \\ 
\Lambda_2 m^{q-p\beta }\left( \int_{B_{R}^c }\left| x\right| ^{\frac{\alpha
-\nu \left( q-p\beta \right) }{1-\beta }}dx\right) ^{1-\beta }\left\|
u\right\| ^{p\beta} \medskip  & \text{if }\frac{1}{p}
<\beta <1 \\ 
\Lambda_2 m^{q-p}\left( \int_{B_{R}^c}\left| x\right| ^{\frac{p}{p-1}(\alpha -\nu \left(
q-p\right)) }V\left( \left| x\right| \right) \left| u\right| ^{p}dx\right) ^{
\frac{p-1}{p}}\left\|u\right\|  & \text{if }\beta =1.
\end{array}
\right. 
\]
\end{lem}

In proving Theorem \ref{THM2}, we will also need the pointwise estimate given by the following lemma.

\begin{lem}\label{lemmainfty}
Assume that there exists $R_{2}>0$ such that 
$$
\essinf_{r>R_{2}}r^{\gamma _{\infty }}V\left(r\right) >0\quad \text{for some }\gamma_{\infty }\leq p -a_{\infty}.
$$

\noindent Then there exists a constant $m_2>0$ such that 
\begin{equation}
\forall u\in X ,\quad \left|
u\left( x\right) \right| \leq 
m_2
\left\| u\right\| \left| x\right| ^{-\frac{p(N-1)-\gamma_{\infty } (p-1)+a_{\infty}}{p^2}%
}\quad \text{almost everywhere in }B_{R_{2}}^{c}.  \label{PointwiseInfty}
\end{equation}
\end{lem}

\proof
Define 
\[
\lambda_2:=
\essinf_{r>R_{2}}r^{\gamma _{\infty }}V\left(r\right)
\]
for brevity. Let $u\in X$ and denote by $C$ any positive constant which does not depend on $u$. 
Let $\tilde u : \left(0, \infty \right)\rightarrow \mathbb{R}^{N}$ be such that $u(x) = \tilde u (|x| )$. Take any $0<a<b<+\infty$, $I=(a,b)$, $I_N = \left\{ x \in \mathbb{R}^{N} \, | \, |x| \in I \right\}$. 
Thanks to $\left( \mathbf{A}\right) $ we know that $|\nabla u| \in L^p (I_N )$, and thanks to the pointwise estimates of Section \ref{SEC:MAIN} we know that $u\in L^{\infty} (I_N )$. Then, using arguments similar to those of Lemma 27 of \cite{BGR_II}, it is easy to see that $\tilde u \in W^{1,p} (I )$. Define now

$$v(r):= r^{\frac{(N-1)p- \gamma_{\infty}(p-1) + a_{\infty}}{p}} \, |\tilde u (r) |^p ,$$

\noindent and $\lambda := \liminf_{r\rightarrow +\infty} v(r)$. Of course $\lambda \geq 0$. As a first step, we want to prove that $\lambda =0$, and we argue by contradiction. So assume $\lambda >0$. Hence there exists ${\overline r} >R_2 $ such that for all $r \geq {\overline r}$ it holds $v(r) \geq \lambda/2$ that is

$$r^{\frac{(N-1)p- \gamma_{\infty}(p-1) + a_{\infty}}{p}} \, |\tilde u (r) |^p \geq \frac{\lambda}{2}.$$

\noindent Hence

\begin{eqnarray*}
\int_{B_{R_2}^c} V(|x|) \, |u(x)|^p dx &\geq &\lambda_{2} \int_{B_{R_2}^c}\frac{|u(x)|^p}{|x|^{\gamma_{\infty}}} \, dx 
= C 
\int_{R_2}^{+\infty}\frac{|\tilde u (r)|^p}{r^{\gamma_{\infty}}} \, r^{N-1} \, dr \\
&\geq &\int_{\overline r}^{+\infty}\frac{\lambda}{2 r^{ \frac{\gamma_{\infty}+a_{\infty}}{p}}}dr =+\infty ,
\end{eqnarray*}
because $\gamma_{\infty}+a_{\infty} \leq p$ by hypothesis. This contradiction proves that $\lambda =0$. Hence there exists a sequence $\{ r_n \}_n$ such that 
$r_n \rightarrow +\infty$ and $v(r_n ) \rightarrow 0$. Now let us fix $r>R_2$ and $r_n >r$, and consider the interval $I=(r, r_n )$. We know that ${\tilde u } \in
D^{1,p}(I) \cap L^{\infty} (I)$, from which we get that $|{\tilde u}|^p $ has a weak derivative in $I$ given by
$$\frac{d}{dr} |{\tilde u}(r)|^p = p | {\tilde u} (r) |^{p-2}\, {\tilde u}(r) \,{\tilde u}'(r) .$$

\noindent Furthermore, the function $r \rightarrow r^{\frac{(N-1)p- \gamma_{\infty}(p-1) + a_{\infty}}{p}}$ is a bounded smooth function in $I$ with bounded derivative, from which we derive that the weak derivatives of $v$ can be computed with the usual Leibniz rule, as follows

$$v'(s) =\frac{(N-1)p- \gamma_{\infty}(p-1) + a_{\infty}}{p}\,  s^{\left( N-2 - \gamma_{\infty} \frac{p-1}{p}   + \frac{a_{\infty}}{p}\right) }\, |{\tilde u} (s)|^p $$

$$+
p \, s^{\frac{(N-1)p- \gamma_{\infty}(p-1) + a_{\infty}}{p}} \, | {\tilde u} (s) |^{p-2}\, {\tilde u}(s) \,{\tilde u}'(s).
$$

\noindent From this it is easy to derive that $v' \in L^p (I)$, so that $v \in W^{1,p}(I)$. From the hypotheses $\gamma_{\infty} \leq p-a_{\infty}$ and $a_{\infty} >p-N$, we deduce that $\gamma_{\infty} (p-1) -a_{\infty} \leq (N-1)p$ and hence

$$ v'(s) \geq p \, s^{\frac{(N-1)p- \gamma_{\infty}(p-1) + a_{\infty}}{p}} | {\tilde u} (s) |^{p-2}\, {\tilde u}(s) \,{\tilde u}'(s) \geq 
- p \, s^{\frac{(N-1)p- \gamma_{\infty}(p-1) + a_{\infty}}{p}} | {\tilde u} (s) |^{p-1}\, |{\tilde u}'(s) |. $$

\noindent As $v\in W^{1,p}(I)$, we can apply the fundamental theorem of Calculus to get

$$v(r_n ) -v(r) = \int_{r}^{r_n} v' (s) ds \geq -p\int_{r}^{r_n}  s^{\frac{(N-1)p- \gamma_{\infty}(p-1) + a_{\infty}}{p}} | {\tilde u} (s) |^{p-1}\,
 |{\tilde u}'(s) |ds, $$

\noindent whence

$$v(r ) -v(r_n ) \leq p\int_{r}^{r_n}  s^{\frac{(N-1)p- \gamma_{\infty}(p-1) + a_{\infty}}{p}} | {\tilde u} (s) |^{p-1}\,
 |{\tilde u}'(s) |ds $$

$$ =p \int_{r}^{r_n} s^{(N-1) \, \frac{p-1}{p}} \frac{|{\tilde u}(s)|^{p-1}}{s^{\gamma_{\infty} \frac{p-1}{p}}} \,  s^{\frac{N-1}{p}} \, 
s^{\frac{a_{\infty}}{p}} \, |u'(s)| \, ds 
$$

$$
\leq p\left(  \int_{r}^{r_n} s^{N-1}  \frac{|{\tilde u}(s)|^{p}}{s^{\gamma_{\infty} }} \, ds  \right)^{\frac{p-1}{p}} \, \left( \int_{r}^{r_n} s^{N-1} s^{a_{\infty}}\, |u'(s)|^p \, ds      \right)^{1/p} 
$$

$$
\leq p\left(  \int_{R_2}^{+\infty} s^{N-1}  \frac{1}{\lambda_{2} } V(s)|{\tilde u}(s)|^{p} \, ds  \right)^{\frac{p-1}{p}} \, \left(C \int_{R_2}^{+\infty} s^{N-1}
 A(s)\, |u'(s)|^p \, ds      \right)^{1/p} 
$$

$$
= C 
\left( \int_{B_{R_2}^c} V(|x|) |u(x)|^p dx    \right)^{\frac{p-1}{p}} \, 
\left( \int_{B_{R_2}^c} A(|x|) |\nabla u(x)|^p dx    \right)^{\frac{1}{p}} 
\leq C ||u||^p .$$

%
%
%

\noindent As $C$ does not depend on $n$, we can pass to the limit as $n\rightarrow \infty$ and get $v(r) \leq C ||u||^p$, that is 

$$r^{\frac{(N-1)p- \gamma_{\infty}(p-1) + a_{\infty}}{p}} \, |\tilde u (r) |^p \leq C ||u||^p .$$

\noindent This gives the thesis. 

\endproof

For future reference, we define three functions $\alpha
_{1}\left( \beta ,\gamma \right) $, $\alpha _{2}\left( \beta
\right) $ and $\alpha _{3}\left( \beta ,\gamma \right) $ by
setting 
\begin{equation}
\alpha _{1}\left( \beta ,\gamma \right):=-\left( 1-\beta \right) \gamma ,\quad \alpha _{2}\left( \beta
\right):=-\left(
1-\beta \right) N,\quad \alpha _{3}\left( \beta ,\gamma \right):=-\frac{(p-1) N+\left( 1-p\beta \right) \gamma 
}{p}.  \label{alpha_i :=}
\end{equation}

\proof[Proof of Theorem \ref{THM2}]
Assume the hypotheses of the theorem and denote 
\[
\Lambda _{2}:=\esssup_{r>R_{2}}\frac{K\left( r\right) }{%
r^{\alpha _{\infty }}V\left( r\right) ^{\beta _{\infty }}}\quad \text{and}%
\quad \lambda _{2}:=\essinf_{r>R_{2}}r^{\gamma _{\infty
}}V\left( r\right) . 
\]
Observe that $\forall \xi \geq 0$ one has 
\begin{equation}
\esssup_{r>R_2}\frac{K\left( r\right) }{r^{\alpha _{\infty }+\xi
\gamma _{\infty }}V\left( r\right) ^{\beta _{\infty }+\xi }}
= \esssup_{r>R_{2}}\frac{K\left( r\right) }{r^{\alpha _{\infty }}V\left(
r\right) ^{\beta _{\infty }}\left( r^{\gamma _{\infty }}V\left( r\right)
\right) ^{\xi }}\leq \frac{\Lambda _{2}}{\lambda _{2}^{\xi }}
<+\infty .  \label{stimaETA}
\end{equation}
Let $u\in X $ be such that $\left\|u\right\| =1$ and let $R\geq R_{2}$. 
We will denote by $C$ any positive constant which does not depend on $u$
or $R$ (such as $\Lambda _{2}/\lambda _{2}^{\xi }$ if $\xi $
does not depend on $u$ or $R$).
We will distinguish several cases, in each of which
we will choose a suitable $\xi \geq 0$ and, thanks to (\ref{stimaETA}) and (\ref{PointwiseInfty}), we will apply Lemma 
\ref{Lem(Omega2)} with $\alpha =\alpha _{\infty }+\xi \gamma _{\infty }$, $\beta
=\beta _{\infty }+\xi $ (whence $\Lambda_2 $ will be given by the left hand
side of (\ref{stimaETA})), $m=m_2\left\| u\right\| =m_2$ and 
$\nu =\frac{p(N-1)-\gamma _{\infty }(p-1) +a_{\infty}}{p^2}$. 
Recall that we are assuming $\gamma _{\infty }\leq p - a_{\infty}$ and $a_{\infty} >p -N$, so that $\nu>0$. 
We will obtain that 
\[
\int_{B_{R}^{c}}K\left( \left| x\right| \right) \left| u\right|
^{q_{2}} dx\leq CR^{\delta }
\]
for some $\delta <0$, not dependent on $R$, so that the result
follows.
\par \noindent We set $\alpha_1 =\alpha
_{1}\left( \beta_{\infty} ,\gamma_{\infty} \right) $, $\alpha_2 =\alpha _{2}\left( \beta_{\infty}
\right) $ and $\alpha_3 =\alpha _{3}\left( \beta_{\infty} ,\gamma_{\infty} \right) $.

\medskip 

\noindent \emph{Case }$\alpha _{\infty }\geq \alpha _{1}$.\smallskip

\noindent We take $\xi =1-\beta _{\infty }$ and apply Lemma \ref{Lem(Omega2)}
with $\beta =\beta _{\infty }+\xi =1$ and $\alpha =\alpha _{\infty }+\xi
\gamma _{\infty }=\alpha _{\infty }+\left( 1-\beta _{\infty }\right) \gamma
_{\infty }$. We get 
\begin{eqnarray*}
\int_{B_{R}^{c}}K\left( \left| x\right| \right) \left| u\right|
^{q_{2}} dx &\leq &C\left( \int_{B_{R}^{c}}\left| x\right|
^{\frac{p}{p-1}\left( \alpha -\nu \left( q_{2}-p\right) \right) }V\left( \left| x\right| \right)
\left| u\right| ^{p}dx\right) ^{\frac{p-1}{p}} \\
&\leq &C\left( R^{\frac{p}{p-1} \left( \alpha -\nu \left( q_{2}- p \right) \right) 
}\int_{B_{R}^{c}}V\left( \left| x\right| \right) \left| u\right|
^{p}dx\right) ^{\frac{p-1}{p}}\leq CR^{\alpha -\nu \left( q_{2}-p\right) }
\end{eqnarray*}
where we used the fact that $\alpha -\nu \left( q_{2}-p\right) <0$. Indeed
$$
\alpha - \nu (q_2 -p) = \alpha_{\infty} + \gamma_{\infty} (1- \beta_{\infty}) -\frac{p(N-1)-\gamma _{\infty }(p-1) +a_{\infty}}{p^2} (q_2 -p) 
$$
$$
= \frac{p(N-1)-\gamma _{\infty }(p-1) +a_{\infty}}{p^2} \left(  \frac{p^2 (\alpha_{\infty}+\gamma_{\infty} (1- \beta_{\infty})) }{p(N-1)-\gamma _{\infty }(p-1) +a_{\infty}}  -q_2 +p     \right),
$$
\noindent where
$$
 \frac{p^2 (\alpha_{\infty}+\gamma_{\infty} (1- \beta_{\infty})) }{p(N-1)-\gamma _{\infty }(p-1) +a_{\infty}} \, +p $$
$$=\frac{p}{p(N-1)-\gamma _{\infty }(p-1) +a_{\infty}}
\left(  p\alpha_{\infty} + p  \gamma_{\infty} - p \gamma_{\infty}  \beta_{\infty} + p (N-1) - p  \gamma_{\infty}  +   \gamma_{\infty} + a_{\infty}    \right)$$
$$
=q_{**} (\alpha_{\infty} , \beta_{\infty} , \gamma_{\infty} ,a_{\infty} ),
$$
\noindent so that 
$$
\alpha - \nu (q_2 -p) =\frac{p(N-1)-\gamma _{\infty }(p-1) +a_{\infty}}{p^2} \left(  q_{**} - q_2 \right) <0,
$$
\noindent because $q_2 > q_{**}$ and $p(N-1)-\gamma _{\infty }(p-1) +a_{\infty}>0$ by hypothesis.

\bigskip

\noindent \emph{Case }$\max \left\{ \alpha _{2},\alpha _{3}\right\} <\alpha
_{\infty }<\alpha _{1}$.\smallskip

\noindent Take $\xi =\frac{\alpha _{\infty }+\left( 1-\beta _{\infty
}\right) N}{N-\gamma _{\infty }}>0$ and apply Lemma \ref{Lem(Omega2)} with $%
\beta =\beta _{\infty }+\xi $ and $\alpha =\alpha _{\infty }+\xi \gamma
_{\infty }$. For doing this, observe that $\alpha _{3}<\alpha _{\infty
}<\alpha _{1}$ implies 
\[
\beta =\beta _{\infty }+\xi =\frac{\alpha _{\infty }-\gamma _{\infty }\beta
_{\infty }+N}{N-\gamma _{\infty }}= q_{*}\left( \alpha_{\infty}, \beta_{\infty},  \gamma_{\infty}      \right) \in \left( \frac{1}{p},1\right) . 
\]
\noindent Notice that $\alpha +N(1-\beta ) =0$. Indeed 
$$
1-\beta = 1- \frac{\alpha_{\infty}   -\beta_{\infty} \gamma_{\infty} +N}{N-\gamma_{\infty}}= \frac{\beta_{\infty} \gamma_{\infty} -\gamma_{\infty} - \alpha_{\infty}}{N-\gamma_{\infty}}
$$
\noindent while
$$\alpha = \alpha_{\infty} + \xi \gamma_{\infty} =\alpha_{\infty} +\gamma_{\infty}\frac{\alpha _{\infty }+\left( 1-\beta _{\infty
}\right) N}{N-\gamma _{\infty }}= \frac{N \alpha_{\infty} - \alpha_{\infty} \gamma_{\infty} + \alpha_{\infty} \gamma_{\infty} +N \gamma_{\infty} -N \gamma_{\infty} \beta_{\infty}}{N-\gamma_{\infty}} $$
$$=\frac{N}{N-\gamma_{\infty}} \left(  -\beta_{\infty} \gamma_{\infty} + \gamma_{\infty} + \alpha_{\infty}    \right),
$$
\noindent so that $\alpha +N(1-\beta ) =0$ easily follows. Now from Lemma \ref{Lem(Omega2)} we get  
\[
\int_{B_{R}^{c}}K\left( \left| x\right| \right) \left| u\right|
^{q_{2}} dx\leq C\left( \int_{B_{R}^{c}}\left| x\right| ^{%
\frac{\alpha -\nu \left( q_{2}-p\beta \right) }{1-\beta }}dx\right)
^{1-\beta }\leq C\left( R^{\frac{\alpha -\nu \left( q_{2}-p\beta \right) }{%
1-\beta }+N}\right) ^{1-\beta },
\]
\noindent because 
$$
\frac{\alpha -\nu \left( q_{2}-p\beta \right) }{1-\beta }+N =\frac{1}{1-\beta}
\left(  \alpha -\nu (q_2 -p\beta ) +N(1-\beta )     \right) =\frac{\nu}{1-\beta} \left(  p \beta - q_2 \right)$$
$$=\frac{\nu }{1-\beta }\left(
q_{*}(\alpha_{\infty}, \beta_{\infty}, \gamma_{\infty})-q_{2}\right) <0$$
\noindent by hypothesis $q_2 > q_{*}$. 
\bigskip

\noindent \emph{Case }$\beta _{\infty }=1$\emph{\ and }$\alpha _{\infty
}\leq 0=\alpha _{2}\,\left( =\max \left\{ \alpha _{2},\alpha _{3}\right\}
\right) $.\smallskip

\noindent Take $\xi =0$ and apply Lemma \ref{Lem(Omega2)} with $\beta =\beta
_{\infty }+\xi =1$ and $\alpha =\alpha _{\infty }+\xi \gamma _{\infty
}=\alpha _{\infty }$. We get 
\[
\int_{B_{R}^{c}}K\left( \left| x\right| \right) \left| u\right|
^{q_{2}} dx\leq C\left( \int_{B_{R}^{c}}\left| x\right|
^{\frac{p}{p-1}\left( \alpha _{\infty }-\nu \left( q_{2}-p\right) \right) }V\left( \left| x\right|
\right) \left| u\right| ^{p}dx\right) ^{\frac{p-1}{p}}\leq CR^{\alpha _{\infty
}-\nu \left( q_{2}-p\right) }, 
\]
since $\alpha _{\infty }-\nu \left( q_{2}-p\right) \leq -\nu \left(
q_{2}-p\right) <0.$
\bigskip

\noindent \emph{Case }$\frac{1}{p}<\beta _{\infty }<1$\emph{\ and }$\alpha
_{\infty }\leq \alpha _{2}\,\left( =\max \left\{ \alpha _{2},\alpha
_{3}\right\} \right) $.\smallskip

\noindent Take $\xi =0$ again and apply Lemma \ref{Lem(Omega2)} with $\beta
=\beta _{\infty }\in \left( \frac{1}{p},1\right) $ and $\alpha =\alpha
_{\infty }$. We get 
\[
\int_{B_{R}^{c}}K\left( \left| x\right| \right) \left| u\right|^{q_{2}}dx 
\leq  C\left( \int_{B_{R}^{c}}\left| x\right| ^{\frac{\alpha _{\infty }-\nu \left( q_{2}-p\beta _{\infty }\right) }{1-\beta
_{\infty }}}dx\right) ^{1-\beta _{\infty }}
\leq  C\left( R^{\frac{\alpha_{\infty }-\nu \left( q_{2}-p\beta _{\infty }\right) }{1-\beta _{\infty }}%
+N}\right) ^{1-\beta _{\infty }}, 
\]
since 
\[
\frac{\alpha _{\infty }-\nu \left( q_{2}-p\beta _{\infty }\right) }{1-\beta
_{\infty }}+N=\frac{\alpha _{\infty }-\alpha _{2}-\nu \left( q_{2}-p\beta
_{\infty }\right) }{1-\beta _{\infty }}<0 .
\]
\smallskip

\noindent \emph{Case }$\beta _{\infty }\leq \frac{1}{p}$\emph{\ and }$\alpha
_{\infty }\leq \alpha _{3}\,\left( =\max \left\{ \alpha _{2},\alpha
_{3}\right\} \right) $.\smallskip

\noindent Take $\xi =\frac{1-p\beta _{\infty }}{p}\geq 0$, we can apply
Lemma \ref{Lem(Omega2)} with $\beta =\beta _{\infty }+\xi =\frac{1}{p}$ and $%
\alpha =\alpha _{\infty }+\xi \gamma _{\infty }$. We get 
\[
\int_{B_{R}^{c}}K\left( \left| x\right| \right) \left| u\right|
^{q_{2}}dx
\leq 
C\left( \int_{B_{R}^{c}}
\left| x\right|^{\frac{ \left( \alpha -\nu (q_2-1)\right)p}{p-1} }    dx\right)
^{\frac{p-1}{p}}
\leq
CR^{\alpha -\nu \left( q_{2}-1\right) +\frac{N(p-1)}{p}}, 
\]
since 
\[
\alpha -\nu \left( q_{2}-1\right) +\frac{N(p-1)}{p}=\alpha _{\infty }+\frac{%
1-p\beta _{\infty }}{p}\gamma _{\infty }+\frac{N(p-1)}{p}-\nu \left(
q_{2}-1\right) =\alpha _{\infty }-\alpha _{3}-\nu \left( q_{2}-1\right) <0. 
\]

\endproof

The proof of Theorem \ref{THM3} will be achieved by the following further lemmas.

\begin{lem}
\label{LEM(pointwise0)}Assume that there exists $R_1>0$ such that 
\[
\essinf_{r\in \left( 0,R_1\right)
}r^{\gamma _{0}}V\left( r\right) >0\quad \text{for some }\gamma _{0}\geq p -a_0.
\]
Then there exists a constant $m_1 > 0$ 
such that $\forall u\in X$ 
one has 
\begin{equation}
\left| u\left( x\right) \right| \leq 
m_1 \left\| u\right\| \left| x\right|^{-\frac{
p(N-1)-(p-1)\gamma _{0}+a_{0}}{p^2}} \quad \text{almost everywhere in }B_{R_1}.
\label{Pointwise0}
\end{equation}
\end{lem}

\proof%
We argue as in Lemma \ref{lemmainfty}. Define 
\[
\lambda_1:=
\essinf_{r\in \left( 0,R_1\right)}r^{\gamma _{0}}V\left( r\right)
\quad\text{and}\quad
\mu_1  :=\sup_{r\in \left( 0,R_1\right)}   \frac{  r^{a_{0}}   }{   A\left( r\right)   }
\]
for brevity. Let $u\in X$ and denote by $C$ any positive constant which does not depend on $u$. 
Let ${\tilde u}:(0, +\infty ) \rightarrow \mathbb R$ be such that $u(x) = {\tilde u} (|x|)$. As before, for any $0<a<b<+\infty$, if 
$I$ and $I_N$ are as in Lemma \ref{lemmainfty} we obtain $u \in W^{1,p}(I_N ) $ and ${\tilde u} \in W^{1,p}(I ) $. Define now
$$ v(r) = r^{N-1 - \frac{p-1}{p}\gamma_0 + \frac{a_0}{p}} \, |{\tilde u} (r)|^p  , \quad \quad l_0 = \liminf_{r\rightarrow 0^+} v(r) .$$
\noindent Of course $l_0 \geq 0$. We want to prove $l_0 =0$, and we argue by contradiction. So assume $l_0 >0$. Then there exists ${\overline r } <R_1$ such that, for all $0<r< {\overline r }$, one has
$$ |{\tilde u} (r)|^p  \, r^{N-1 - \gamma_0 } \geq \frac{l_0}{2} \frac{1}{r^{\frac{\gamma_0 + a_0}{p}}}$$
\noindent and hence
$$\int_{B_{R_1}} V(|x|) |u(x)|^p dx \geq \lambda_1 \int_{B_{R_1}} \frac{|u(x)|^p}{|x|^{\gamma_0}}dx = 
C\int_{0}^{R_!} |{\tilde u} (r)|^p \, r^{N-1 -\gamma_0} dr  $$
$$\geq C 
\frac{l_0}{2} \int_{0}^{{\overline r}} \frac{1}{r^{\frac{\gamma_0 + a_0}{p}}} dr = +\infty ,$$
\noindent because $\gamma_0 + a_0 \geq p$ by hypothesis. This contradiction proves that $l_0 =0$. Hence, there exists a sequence $\{ r_n \}_n $ such that $r_n \rightarrow 0$ and $v(r_n ) \rightarrow 0$. Now fix $r<R_1$ and $r_n <r$, and set $I=(r_n , r)$. Arguing as in Lemma \ref{lemmainfty}, we have that $v \in W^{1,p}(I)$ and the derivative of $v$ is given by
$$
v'(s) = \left( N-1 - \frac{p-1}{p} \gamma_0 + \frac{a_0}{p}   \right) \, s^{ N-2 - \frac{p-1}{p} \gamma_0 + \frac{a_0}{p} } \, |{\tilde u} (s)|^p +$$
$$
+p\, s^{ N-1 - \frac{p-1}{p} \gamma_0 + \frac{a_0}{p} } |{\tilde u} (s) |^{p-2}{\tilde u} (s) {\tilde u}' (s) .
$$
\noindent We then apply the fundamental theorem of Calculus to get
$$
v(r)-v(r_n )= \int_{r_n}^{r} v'(s) ds .
$$
\noindent We notice now that 
$$
\int_{r_n}^{r} s^{ N-2 - \frac{p-1}{p} \gamma_0 + \frac{a_0}{p} } \, |{\tilde u} (s)|^p ds =
 \int_{r_n}^{r} s^{N-1} \, s^{\frac{\gamma_0 + a_0 -p}{p} } \, 
\frac{|{\tilde u} (s)|^p }{s^{\gamma_0}} ds \leq r^{\frac{\gamma_0 + a_0 -p}{p} } \int_{0}^{R_1}\frac{|{\tilde u} (s)|^p }{s^{\gamma_0}} s^{N-1} ds 
$$
$$\leq r^{\frac{\gamma_0 + a_0 -p}{p} } \frac{1}{\lambda_1} \int_{0}^{R_!} V(s)|{\tilde u} (s)|^p s^{N-1} ds = 
r^{\frac{\gamma_0 + a_0 -p}{p} } C 
\int_{B_{R_1}} V(|x|) |u(x)|^p dx \leq 
C R_1^{\frac{\gamma_0 + a_0 -p}{p} }||u||^p .$$
\noindent On the other hand
$$
\int_{r_n}^{r} s^{ N-1 - \frac{p-1}{p} \gamma_0 + \frac{a_0}{p} } |{\tilde u} (s) |^{p-2} {\tilde u} (s) {\tilde u}' (s) ds \leq
\int_{r_n}^{r} s^{ (N-1 -  \gamma_0 ))\frac{p-1}{p}} |{\tilde u} (s) |^{p-1} \,s^{ (N-1 + a_0 ))\frac{1}{p}} |{\tilde u}' (s) | ds
$$
$$\leq \left( \int_{r_n}^{r} s^{N-1} \frac{ |{\tilde u} (s) |^{p}}{s^{\gamma_0 }}  ds  \right)^{\frac{p-1}{p}} \, \left( \int_{r_n}^{r} s^{a_0}  |{\tilde u}' (s) |^p s^{N-1} ds \right)^{1/p} 
$$
$$
\leq \left(\frac{1}{\lambda_1} \int_{r_n}^{r}  s^{N-1} V(s)  |{\tilde u} (s) |^{p} ds  \right)^{\frac{p-1}{p}} \, \left(\int_{r_n}^{r} \mu_1 A(s)   |{\tilde u}' (s) |^p s^{N-1} ds \right)^{1/p} 
$$
$$
\leq C 
\left( \int_{B_{R_1}} V(|x|) |u(x)|^p dx   \right)^{\frac{p-1}{p}} \,
\left(   \int_{B_{R_1}} A(|x|) |\nabla u(x)|^p dx   \right)^{1/p} \leq 
C ||u||^p .
$$
\noindent We also notice that from the hypothesis $a_0 + \gamma_0 \geq p$ one deduces $N-1 - \frac{p-1}{p} \gamma_0 + \frac{a_0}{p}  \leq N-p +a_0 $, and hence 
$$v(r)-v(r_n )= \int_{r_n}^{r} v'(s) ds  $$
$$=\left( N-1 - \frac{p-1}{p} \gamma_0 + \frac{a_0}{p}   \right) \int_{r_n}^{r} s^{ N-2 - \frac{p-1}{p} \gamma_0 + \frac{a_0}{p} } \, |{\tilde u} (s)|^p   ds + $$
$$+ p \int_{r_n}^{r} s^{ N-1 - \frac{p-1}{p} \gamma_0 + \frac{a_0}{p} } |{\tilde u} (s) |^{p-2}{\tilde u} (s) {\tilde u}' (s)   ds $$
$$
\leq \left( N-p+a_0 \right) \int_{r_n}^{r} s^{ N-2 - \frac{\gamma_0}{p}  + \frac{a_0}{p} }\,s^{-\gamma_0} \, |{\tilde u} (s)|^p   ds+C\, ||u||^p 
$$
$$\leq
\left( N-p+a_0 \right) R_1^{\frac{\gamma_0 + a_0 -p}{p} } 
\int_{r_n}^{r} s^{ N-1}\,\frac{1}{\lambda_1} \, V(s) \, |{\tilde u} (s)|^p   ds + C\, ||u||^p \leq C\, ||u||^p .
$$
%
%
\noindent This holds for all $n$ large enough (i.e. such that $r_n<r$), so we can pass to the limit as $n\rightarrow\infty$ and get
$v(r)\leq C\, ||u||^p$, which yields the thesis. 
%
%
%
\endproof

\begin{lem}
\label{Lem(Omega0)}Assume that there exists $R_1>0$ such that 
\begin{equation}
\esssup_{r\in \left(0,R_1\right) }\frac{K\left( r\right) }{r^{\alpha }V\left( r\right) ^{\beta }}%
<+\infty \quad \text{for some }\frac{1}{p}\leq \beta \leq 1\text{~and }%
\alpha \in \mathbb{R}  \label{Lem(Omega0): hp1}
\end{equation}
and 
\[
\essinf_{r\in \left( 0,R_1\right)}r^{\gamma _{0}}V\left( r\right) >0\quad \text{for some }\gamma _{0}>p -a_{0}.
\]
Assume also that there exists $q>p\beta $ such that 
\[
\left( p(N-1)-(p-1)\gamma _{0} +a_0 \right) q < p^2 (\alpha +N)-p \beta \left((p-1) \gamma _{0}+p -a_0 \right).
\]
Then there exists a constant $\tilde{C}_1>0$ such that $\forall R\in (0,R_1)$ and $\forall u \in X$ one has 
\[
\int_{B_{R}}K\left( \left| x\right| \right) \left| u\right| ^{q} dx\leq  \tilde{C}_1 
R^{\frac{p^2 (\alpha +N)-p \beta \left( (p-1) \gamma _{0}+p -a_0 \right) -\left( p(N-1)-(p-1)\gamma _{0} +a_0\right) q}{p^2}
}\left\| u\right\| ^{q} .
\]
\end{lem}

\proof
Let $u\in X $ and $R\in (0,r_1)$. 
Denote by $C$ any positive constant which does not depend on $u$ or $R$.
By assumption (\ref{Lem(Omega0): hp1}) and Lemma \ref{LEM(pointwise0)}, we can apply Lemma \ref{Lem(Omega)} with 
$\nu =\frac{(N-1)p-(p-1)\gamma _{0} +a_0}{p^2}$ and $m=m_1 \left\|u\right\|$.
If $\frac{1}{p}\leq \beta <1$, we get 
\begin{eqnarray*}
\int_{B_{R}}K\left( \left| x\right| \right) \left| u\right| ^{q} dx 
&\leq & C m^{q-p\beta }\left( \int_{B_R }\left| x\right|
^{\frac{\alpha -\nu \left( q-p\beta \right) }{1-\beta }}dx  \right) ^{1-\beta
}\left\| u\right\| ^{p\beta } 
= C 
\left( \int_{B_R }\left| x\right|
^{\frac{\alpha -\nu \left( q-p\beta \right) }{1-\beta }}dx  \right) ^{1-\beta }\left\| u\right\|^{q}. 
\end{eqnarray*}
Notice now that
$$ \frac{\alpha - \nu (q-p\beta)}{1-\beta}= \frac{1}{1-\beta} \left( \alpha - \frac{p(N-1) -\gamma_0 (p-1) +a_0}{p^2} \left( q-p\beta \right)     \right) 
$$
$$
=\frac{\alpha p^2 +p\beta \left( p(N-1) -\gamma_0 (p-1) +a_0  \right) -q \left( p(N-1) -\gamma_0 (p-1) +a_0  \right) }{p^2 (1-\beta)}, 
$$
\noindent so that
$$ \frac{\alpha - \nu (q-p\beta)}{1-\beta} +N = \frac{p^2 ( \alpha +N) -p\beta \left( \gamma_0 (p-1) +p -a_0  \right) -q \left( p(N-1) -\gamma_0 (p-1) +a_0  \right) }{p^2 (1-\beta )} >0
$$
\noindent by hypothesis. Hence we get 
$$
\int_{B_R}\left| x\right|
^{\frac{\alpha -\nu \left( q-p\beta \right) }{1-\beta }}dx =  C R^{\frac{ p^2 ( \alpha +N) -p\beta \left( \gamma_0 (p-1) +p -a_0  \right) -q \left( p(N-1) -\gamma_0 (p-1) +a_0  \right)}{p^2 (1-\beta )} }, 
$$
\noindent and from this we get the thesis. 
If instead we have $\beta =1$, we get 
\begin{eqnarray*}
&&\int_{B_{R}}K\left( \left| x\right| \right) \left| u\right| ^{q} dx \\
&\leq &\Lambda m^{q-p}\left( \int_{B_R }\left| x\right|
^{\frac{p}{p-1}\left( \alpha -\nu \left( q-p\right)\right) }V\left( \left| x\right| \right) \left|
u\right| ^{p}dx\right) ^{\frac{p-1}{p}}\left\| u\right\| \\
&=& C \left( \int_{B_{R}}\left| x\right| ^{\frac{
p^2 \alpha -\left( p(N-1)- (p-1)\gamma _{0} +a_0 \right) \left( q-p\right) }{p(p-1)}}V\left(
\left| x\right| \right) \left| u\right| ^{p}dx\right) ^{\frac{p-1}{p}}\left\|
u\right\| ^{q-p}\left\| u\right\| \\
&\leq & C \left( R^{\frac{p^2 (\alpha +N) - p \left(
p+(p-1)\gamma _{0} -a_0 \right) -q \left( p(N-1) -\gamma_0 (p-1) +a_0 \right) }{p(p-1)}}  \int_{B_{R}}V\left( \left|
x\right| \right) \left| u\right| ^{p}dx\right) ^{\frac{p-1}{p}}\left\|
u\right\| ^{q-p}\left\| u\right\| \\
&\leq & C R^{\frac{p^2 (\alpha +N) - p \left(
p+(p-1)\gamma _{0} -a_0 \right) -q \left( p(N-1) -\gamma_0 (p-1) +a_0 \right) }{p^2}}   \left\| u\right\| ^{q},
\end{eqnarray*}
since 
$$\alpha - \nu (q-p)= \alpha - \frac{
(N-1)p-(p-1)\gamma _{0} +a_0}{p^2} (q-p) 
$$
$$ 
=\frac{p^2 (\alpha +N) - p \left(
p+(p-1)\gamma _{0} -a_0 \right) -q \left( p(N-1) -\gamma_0 (p-1) +a_0 \right)}{p^2}>0
$$
\noindent from the hypotheses.
\endproof

\proof[Proof of Theorem \ref{THM3}]
Assume the hypotheses of the theorem and denote 
\[
\Lambda _{0}:=\esssup_{r\in \left( 0,R_{1}\right) }\frac{K\left(
r\right) }{r^{\alpha _{0}}V\left( r\right) ^{\beta _{0}}}\quad \text{and}%
\quad \lambda _{0}:=\essinf_{r\in \left( 0,R_{1}\right) }r^{\gamma
_{0}}V\left( r\right).
\]
If $\gamma _{0}=p -a_0$ the thesis of the theorem is true by Theorem 3.2 of \cite{AVK_I}. Indeed in this case $(\alpha_0 , q_1) \in  \mathcal{A}_{a_0 ,\beta _{0},\gamma _{0}}$ means 
$$\max \{ 1, p\beta_0 \} <q_1 < \min \{ q_{*}, q_{**}  \}$$
\noindent and it is easy to verify that 
$$
q_{*} =q_{**}= p \frac{\alpha_0 -\beta_0 (p-a_0 ) +N}{N-p +a_0}  
$$
\noindent This value coincides with $q^{*}(a_0 , \alpha_0 , \beta_0 )$ of Theorem 3.2 of \cite{AVK_I}, so we have 
$$
q^{*}(a_0 , \alpha_0 , \beta_0 )> \max \{ 1, p\beta_0 \} 
$$
\noindent and Theorem 3.2 of \cite{AVK_I} gives the result. Hence, we can
assume $\gamma _{0}>p-a_0 $ without restriction. 
We claim that for every $0<R\leq R_{1}$ there exists $b\left(
R\right) >0$ such that $b\left( R\right) \rightarrow 0$ as $R\rightarrow
0^{+}$ and 
\[
\int_{B_{R}}K\left( \left| x\right| \right) \left| u\right| ^{q_{1}} dx
\leq b\left( R\right) \left\| u\right\| ^{q_{1}},
\quad 
\forall u \in X,
\]
which clearly gives the result. In order to prove this claim, let $0<R\leq R_{1}$.
Then one has
\begin{equation}
\lambda \left( R\right) :=\essinf_{r\in \left( 0,R\right)
}r^{\gamma _{0}}V\left( r\right) \geq \lambda _{0}>0  \label{stimaETA1}
\end{equation}
and for every $\xi \geq 0$ we have 
\begin{eqnarray}
\esssup_{r\in \left( 0,R\right) }\frac{K\left( r\right) }{r^{\alpha
_{0}+\xi \gamma _{0}}V\left( r\right) ^{\beta _{0}+\xi }}= 
\esssup_{r\in \left( 0,R_{1}\right) }\frac{K\left( r\right) }{r^{\alpha
_{0}}V\left( r\right) ^{\beta _{0}}\left( r^{\gamma _{0}}V\left( r\right)
\right) ^{\xi }}
\leq \frac{\Lambda _{0}}{\lambda _{0}^{\xi }}<+\infty .
\label{stimaETA0}
\end{eqnarray}
Denoting $\alpha _{1}=\alpha _{1}\left( \beta _{0},\gamma _{0}\right) $, $%
\alpha _{2}=\alpha _{2}\left( \beta _{0}\right) $ and $\alpha _{3}=\alpha
_{3}\left( \beta _{0},\gamma _{0}\right) $ as defined in (\ref{alpha_i :=}%
), we will now distinguish five cases, which reflect the five definitions (%
\ref{A:=}) of the set $\mathcal{A}_{a_0 ,\beta _{0},\gamma _{0}}$. For the sake
of clarity, some computations will be displaced in the Appendix.\medskip

\noindent \emph{Case }$p - a_0 <\gamma _{0}<N$\emph{.}\smallskip

\noindent In this case, $\left( \alpha _{0},q_{1}\right) \in \mathcal{A}_{\beta _{0},\gamma _{0}}$ means

$$\max \left\{ 1,p\beta _{0}\right\} <q_{1}<\displaystyle\min \left\{ p\frac{\alpha
_{0}-\beta _{0}\gamma _{0}+N}{N-\gamma _{0}}, p\frac{p\alpha _{0}+\left(
1-p\beta _{0}\right) \gamma _{0}+p(N-1)+a_0}{p(N-1)-(p-1)\gamma _{0}+a_0}\right\}. $$

\noindent Notice that from the hypotheses we have $N-\gamma_0 >0$ and also $\frac{p(N-1) +a_0}{p-1}  >N>\gamma_0$, so that $p(N-1) 
-\gamma_0 (p-1) +a_0 >0$. The inequalities above imply that we can fix $\xi \geq 0$, independent of $R$, 
in such a way that $\alpha =\alpha _{0}+\xi \gamma _{0}$
and $\beta =\beta _{0}+\xi $ satisfy 
\begin{equation}
\frac{1}{p}\leq \beta \leq 1\quad \text{and}\quad p\beta <q_{1}<\frac{
p^2 (\alpha +N)-p\beta \left((p-1) \gamma _{0}+p - a_0\right) }{p(N-1)-(p-1)\gamma _{0} +a_0}
\label{Lem(zero): cond1}
\end{equation}
(see Appendix). Hence, by (\ref{stimaETA0}) and (\ref{stimaETA1}), we can
apply Lemma \ref{Lem(Omega0)} (with $\alpha=\alpha_0+\xi\gamma_0$, $\beta=\beta_0+\xi$ and $q=q_{1}$), so that $\forall u \in X$ we get 
\[
\int_{B_{R}}K\left( \left| x\right| \right) \left| u\right| ^{q_{1}} dx
\leq \tilde{C}_1 R^{\frac{p^2 (\alpha
+N)-p \beta\left( (p-1) \gamma _{0}+p -a_0 \right)  -\left( p(N-1)-(p-1)\gamma _{0}+a_0 \right) q_{1}%
}{p^2}}\left\| u\right\| ^{q_{1}} . 
\]
This gives the result, since $R^{p^2(\alpha +N)-p\left( (p-1) \gamma _{0}+p\right)
\beta -\left( p(N-1)- (p-1)\gamma _{0}\right) q_{1}}\rightarrow 0$ as $R\rightarrow
0^{+}$. 
\medskip

\noindent \emph{Case }$\gamma _{0}=N$. \smallskip

\noindent In this case, $\left( \alpha _{0},q_{1}\right) \in \mathcal{A}
_{a_0 ,\beta _{0},\gamma _{0}}$ means 
\[
\quad \max \left\{ 1,p\beta _{0}\right\} <q_{1}<p\frac{p\alpha
_{0}+\left( 1-p\beta _{0}\right) \gamma _{0}+p(N-1)+a_0 }{p(N-1)-(p-1)\gamma _{0}+a_0 } 
\]
and these conditions still ensure that we can fix $\xi \geq 0$ in such a way
that $\alpha =\alpha _{0}+\xi \gamma _{0}$ and $\beta =\beta _{0}+\xi $
satisfy (\ref{Lem(zero): cond1}) (see Appendix), so that the result ensues
again by Lemma \ref{Lem(Omega0)}.\medskip

\noindent \emph{Case }$N<\gamma _{0}<\frac{p(N-1) +a_0}{p-1}$.\smallskip

\noindent In this case, $\left( \alpha _{0},q_{1}\right) \in \mathcal{A}
_{a_0 \beta _{0},\gamma _{0}}$ means 
\[
\begin{tabular}{l}
$\alpha _{0}>\alpha _{1}$\quad and\smallskip\\
$\displaystyle\max \left\{ 1,p\beta _{0},p\frac{\alpha _{0}-\beta _{0}\gamma _{0}+N}{N-\gamma _{0}}\right\} <q_{1}
<\displaystyle p\frac{p\alpha _{0}+\left( 1-p\beta _{0}\right) \gamma _{0}+p(N-1)+a_0 }{p(N-1)-(p-1)\gamma_{0}+a_0}$
\end{tabular}
\]
and the conclusion then follows as in the former cases (see
Appendix).\medskip

\noindent \emph{Case }$\gamma _{0}=\frac{p(N-1)+a_0}{p-1}$.\smallskip

\noindent In this case, $\left( \alpha _{0},q_{1}\right) \in \mathcal{A}
_{a_0 , \beta _{0},\gamma _{0}}$ means 
\[
\alpha _{0}>\alpha _{1}\quad \text{and}\quad \max \left\{ 1,p\beta _{0},p
\frac{\alpha _{0}-\beta _{0}\gamma _{0}+N}{N-\gamma _{0}}\right\} <q_{1} 
\]
and these conditions ensure that we can fix $\xi \geq 0$ in such a way that $
\alpha =\alpha _{0}+\xi \gamma _{0}$ and $\beta =\beta _{0}+\xi $ satisfy 
\[
\frac{1}{p}\leq \beta \leq 1,\quad q_{1}>p\beta \quad \text{and}\quad
0<p^2 \left( \alpha +N \right) - p\beta\left( (p-1) \gamma _{0}+p-a_0 \right)  
\]
(see Appendix). The result then follows again from Lemma \ref{Lem(Omega0)}.\medskip

\noindent \emph{Case }$\gamma _{0}>\frac{p(N-1)+a_0}{p-1}$.\smallskip

\noindent In this case, $\left( \alpha _{0},q_{1}\right) \in \mathcal{A}
_{a_0 ,\beta _{0},\gamma _{0}}$ means 
\[
q_{1}>\max \left\{ 1,p\beta _{0},p\frac{\alpha _{0}-\beta _{0}\gamma _{0}+N}{
N-\gamma _{0}},p\frac{p\alpha _{0}+\left( 1-p\beta _{0}\right) \gamma
_{0}+p(N-1)+a_0}{p(N-1)- (p-1)\gamma _{0}+a_0}\right\} 
\]
and this condition ensures that we can fix $\xi \geq 0$ in such a way that $
\alpha =\alpha _{0}+\xi \gamma _{0}$ and $\beta =\beta _{0}+\xi $ satisfy 
\[
\frac{1}{p}\leq \beta \leq 1\quad \text{and}\quad q_{1}>\max \left\{ p\beta
,p\frac{p(\alpha+N)-\beta \left((p-1)\gamma _{0}+p-a_0\right)  }{p(N-1)-(p-1)\gamma_{0}+a_0}\right\} 
\]
(see Appendix). The result still follows from Lemma \ref{Lem(Omega0)}.
\endproof


\section{Existence results \label{SEC: ex}}

Let $N \geq 3$ and $1<p<N$. In this section we deduce our existence results
about radial weak solutions to the equation 
\begin{equation}
-\mathrm{div}\left(A(|x|) \, |\nabla u|^{p-2} \nabla u\right) +V\left( \left| x\right| \right) |u|^{p-2}u= K(|x|) f(u) \quad \text{in }\mathbb{R}^{N},  \label{EQg}
\end{equation}
i.e. functions $u\in X$ such that 
\begin{equation}
\int_{\mathbb{R}^{N}} A(|x|) |\nabla u|^{p-2}\nabla u\cdot \nabla h\,dx+\int_{
\mathbb{R}^{N}}V\left( \left| x\right| \right) |u|^{p-2}uh\,dx=\int_{
\mathbb{R}^{N}}K(|x|) f(u) h\,dx,\quad 
\forall h\in X,  \label{weak solution}
\end{equation}
where $A$, $V$ and $K$ are potentials satisfying 
$\left( \mathbf{A}\right) $, $\left( \mathbf{V}\right) $ and $\left( \mathbf{K}\right) $, 
and $X$ is the Banach spaces defined in Section \ref{SEC:MAIN}. 

As concerns the nonlinearity $f$, we set $F(t)= \int_{0}^{t} f(s)ds $ and we will use the following assumptions:

\begin{itemize}
\item[$\left( f_{0}\right) $]  $f: {\mathbb{R}} \rightarrow {\mathbb{R}}$ is a continuous functions such that $f(0)=0$.

\item[$\left( f_{1}\right) $]  $\exists \theta >p$ such that $\theta
F\left( t\right) \leq f\left( t\right) t$ for all 
$t\geq 0;$

\item[$\left( f_{2}\right) $]  $\exists t_{0}>0$ such that $F\left(t_{0}\right) >0$.

\item[$\left( f_{q_{1},q_{2}}\right) $]  $\exists q_1 ,q_2 >p $ and $ M>0$ such that 
$0 \leq 
f\left( t\right) \leq M\min \left\{t^{q_{1}-1}, t^{q_{2}-1}\right\} $ for all $t\geq 0$.

\end{itemize}

\noindent Observe that, if $q_{1}\neq q_{2}$, the double-power growth
condition $\left( f_{q_{1},q_{2}}\right) $ is more stringent than the more
usual single-power one, since it implies $\sup_{t>0}\,\left| f\left(
t\right) \right| /t^{q-1}<+\infty $ for $q=q_{1}$, $q=q_{2}$ and every $q$
in between. On the other hand, we will never require $q_{1}\neq q_{2}$ in $
\left( f_{q_{1},q_{2}}\right) $, so that our results will also concern
single-power nonlinearities as long as we can take $q_{1}=q_{2}$.
The simplest $f\in C\left( 
\mathbb{R};\mathbb{R}\right) $ such that $\left( f_{q_{1},q_{2}}\right) $ holds is 

$$f\left( t\right) =\min \left\{ \left| t\right| ^{q_{1}-2}t,\left| t\right|
^{q_{2}-2}t\right\} ,$$

\noindent which also ensures $\left( f_{1}\right) $ if $
q_{1},q_{2}>p$ (with $\theta =\min \left\{ q_{1},q_{2}\right\} $). Another model example is 
\[
f\left( t\right) =\frac{\left| t\right| ^{q_{2}-2}t}{1+\left| t\right|
^{q_{2}-q_{1}}}\quad \text{with }1<q_{1}\leq q_{2},
\]
which ensures $\left( f_{1}\right) $ if $q_{1}>p$ (with $\theta =q_{1}$). Note that, in
both these cases, also $\left( f_{2}\right) $ holds true. Moreover, both of these functions $f$
become $f\left( t\right) =\left| t\right| ^{q-2}t$ if $q_{1}=q_{2}=q$.

Now we set
$$I\left( u\right) :=\frac{1}{p} \left\| u\right\| ^{p}-\int_{\mathbb{R}
^{N}}  K(|x|) F\left( u\right) dx $$
$$=\frac{1}{p} \int_{\mathbb{R}^{N}} A(|x|) |\nabla u|^p dx + \frac{1}{p} \int_{\mathbb{R}^{N}} V(|x|) |u|^p dx - 
\int_{\mathbb{R}^{N}}  K(|x|) F\left( u\right) dx.
$$

\noindent From the continuous embedding result of Theorem \ref{THM(cpt)} and the
results of \cite{BPR} about Nemytski\u{\i} operators on the sum of Lebesgue
spaces, we have that $I$ is a $C^{1}$ functional on $X$
provided that there exist $q_{1},q_{2}>1$ such that $\left(
f_{q_{1},q_{2}}\right) $ and $\left( \mathcal{S}_{q_{1},q_{2}}^{\prime
}\right) $ hold. In this case, the Fr\'{e}chet derivative of $I$ at any $
u\in X$ is given by 
\begin{equation}
I^{\prime }\left( u\right) h=\int_{\mathbb{R}^{N}}A(|x|)\left( |\nabla u|^{p-2}\nabla
u\cdot \nabla h+V\left( \left| x\right| \right) |u|^{p-2}uh\right) dx-\int_{
\mathbb{R}^{N}}f\left( \left| x\right| ,u\right) h\,dx 
\label{PROP:diff: I'(u)h=}
\end{equation}
for all $h\in X$, and therefore the critical points of $I:X\rightarrow \mathbb{R}$ satisfy (\ref{weak solution}) 
for all $h\in X$. The following theorem has been proved  in \cite{AVK_I}.

\begin{thm}
\label{THM:ex} Assume hypotheses $(f_0 )$, $(f_1 )$, $(f_2 )$ and that there
exist $q_{1},q_{2}>p$ such that $\left( f_{q_{1},q_{2}}\right) $ and $\left( 
\mathcal{S}_{q_{1},q_{2}}^{\prime \prime }\right) $ hold. Assume also $\left( \mathbf{A}\right) $, 
$\left( \mathbf{V}\right) $, $\left( \mathbf{K}\right) $.
\noindent Then the functional $I:X\rightarrow \mathbb{R}$ has a nonnegative
critical point $u\neq 0$.
\end{thm}


\noindent Such existence result relies on assumption $\left(\mathcal{S}_{q_{1},q_{2}}^{\prime \prime }\right) $, which is quite abstract and difficult to check.
However it can be granted in concrete cases through Theorems \ref{THM2} and \ref{THM3}, 
which ensure condition $\left(\mathcal{S}_{q_{1},q_{2}}^{\prime \prime }\right) $ for suitable ranges of
exponents $q_{1}$ and $q_{2}$ by explicit conditions on the potentials.
Thanks to this, we get the following theorem,
which is the main existence theorem of this paper. 
The notations for $\mathcal{A}_{a_{0},\beta _{0},\gamma _{0}}$, $q_{*}$ and $q_{**}$ are those of Theorems \ref{THM2} and \ref{THM3}.

\begin{thm}
\label{THM:ex2}
Let $N\geq 3$ and $1<p<N$. Assume hypotheses $(f_0 )$, $(f_1 )$, $(f_2 )$, $\left( \mathbf{A}\right) $, 
$\left( \mathbf{V}\right) $, $\left( \mathbf{K}\right) $.
Assume that there exist $R_{1},R_2>0$ such that 
$$
\esssup_{r\in \left( 0,R_{1}\right) }\frac{K\left( r\right) }{
r^{\alpha _{0}}V\left( r\right) ^{\beta _{0}}}<+\infty \quad \text{for some }
0\leq \beta _{0}\leq 1\text{~and }\alpha _{0}\in \mathbb{R},  
$$
$$
\essinf_{r\in \left( 0,R_{1}\right) }r^{\gamma _{0}}V\left(
r\right) >0\quad \text{for some }\gamma _{0}\geq p -a_{0},
$$
$$
\esssup_{r>R_{2}}\frac{K\left( r\right) }{r^{\alpha _{\infty
}}V\left( r\right) ^{\beta _{\infty }}}<+\infty \quad \text{for some }0\leq
\beta _{\infty }\leq 1\text{~and }\alpha _{\infty }\in \mathbb{R},
$$
$$
\essinf_{r>R_{2}}r^{\gamma _{\infty }}V\left( r\right) >0\quad 
\text{for some }\gamma _{\infty } \leq p - a_{\infty}.  
$$

\noindent Assume furthermore that there exist $q_1 , q_2 >p$ such that

\begin{itemize}

\item $\left( \alpha _{0},q_{1}\right) \in \mathcal{A}_{a_{0},\beta _{0},\gamma _{0}}.$

\item $q_{2}>\max \left\{ 1,p\beta _{\infty },q_{*},q_{**}\right\}.$

\item $\left( f_{q_{1},q_{2}}\right) $ holds.

\end{itemize}

\noindent Then equation (\ref{EQg}) has a nonnegative nontrivial solution $u \in X$.

\end{thm}

\bigskip
\noindent The proof of Theorem \ref{THM:ex2} derives from the compactness result of the previous sections.

\section{Examples}\label{SEC:EX}

\noindent In this section we give some examples that might help to understand what is new in our results. In particular, the following examples are not included in the results of \cite{Su-Wang}, which inspired and motivated our work. We recall that in that paper the authors prove some compactness theorems for potentials which behave as a power both at the origin and at infinity, which yield 
existence results for equation (\ref{EQ}) in which $f(u)$ is a power or a sum of powers. 
In all our examples, we consider for simplicity a nonlinearity defined by $f(t)= \min \{ t^{q_1 -1} , t^{q_2 -1} \}$, and we will see how to choose $p<q_1 \leq q_2$ in such a way to get existence results for problem (\ref{EQ}).

\begin{exa}

Assume $N\geq 3$ as always, and $1<p<N-1$. We take $A(r) = 1/r$ and $V,K$ as follows:
\[
V(r) = e^{1/r} \quad \mbox{if} \; r \in (0,1), \quad \quad V(r)= r^{-p-1}\quad \mbox{if} \; r \geq 1, 
\]
\[
K(r) = e^{1/r} \quad \mbox{if} \; r \in (0,1), \quad \quad K(r)= 1\quad \mbox{if} \; r \geq 1. 
\]
The results of \cite{Su-Wang} do not apply because $K$ does not have a power growth near the origin. To apply our results, we first notice that $a_0 = a_{\infty} = -1$, so the hypothesis $a_0 , a_{\infty}\in (p-N, p ]$ is satisfied. In order to apply Theorem \ref{THM:ex2} we set $\alpha_{\infty}= \beta_{\infty}=0$ and $\gamma_{\infty}= p+1= p-a_{\infty}$. We obtain
$$
q_{*} \left(\alpha_{\infty}, \beta_{\infty}, \gamma_{\infty} \right) = p\, \frac{N }{N-\gamma_{\infty} }=p \, \frac{N }{N-p-1 }  , $$
$$
q_{**} \left(a_{\infty}, \alpha_{\infty}, \beta_{\infty}, \gamma_{\infty} \right) = p \, \frac{p+1 + p(N-1) -1}{p(N-1) -(p+1) (p-1) -1}=
p \, \frac{N }{N-p-1 },
$$
\noindent and thus we can pick $q_2 > p \, \frac{N }{N-p-1 }>p$. Then we take $\alpha_{0}= 0$, $\beta_{0}=1$ 
and pick any $\gamma_0 > \frac{p(N-1)-1}{p-1}$, which is possible thanks to the assumptions on $V$. So we can choose $q_1 > \max \{ 1, p\beta_0 , q_{*}, q_{**} \}$, where
$$
q_{*}= q_{*}(\alpha_0 , \beta_0 , \gamma_0 ) = p \, \frac{-\gamma_0 +N}{N-\gamma_0}=p, 
$$
$$q_{**} =q_{**} (a_0 , \alpha_0 , \beta_0 , \gamma_0 )=
p\, \frac{(1-p)\gamma_0 +p(N-1) -1}{p(N-1)-\gamma_0 (p-1) -1}=p,
$$
\noindent which means $q_1 >p$. As a conclusion, we can choose any $q=q_1 =q_2 > p \, \frac{N }{N-p-1 }$, and we get a nonnegative radial solution to the equation
$$
-\mathrm{div}\left( \frac{1}{|x|}|\nabla u|^{p-2} \nabla u\right) +V\left( \left| x\right| \right) u^{p-1}= K(|x|) u^{q-1}.
$$
\end{exa}

\begin{exa}
Assume $N\geq 4$ and $1<p<N-2$. Notice that $p+1 < \frac{p(N-1)-1}{p-1} $ and let $\gamma_0 \in \left( p+1 , \frac{p(N-1)-1}{p-1}  \right) $. We consider $d>1/2$ large enough (see below) and we take $A, V, K$ as follows:
$$
A(r)=  \min \left\{ r^{-2}, r^{-1} \right\} , \quad V(r)= \max \left\{ r^{-\gamma_0} , r^{1/2} \right\} , \quad K(r)= \max \left\{ r^d , r^{1/2}  \right\}  .
$$
\noindent By a careful analysis of the hypotheses therein, one can see that the results of \cite{Su-Wang} do not apply, at least for large $d'$s.
\par \noindent To apply our results, we set $\beta_0 = \beta_{\infty}  =0$, $\alpha_0 = 1/2$, $\alpha_{\infty} = d$, $a_0 =-1$, $a_{\infty}=-2$ and $\gamma_{\infty}=-1/2$, and let $\gamma_0$ as before. Notice that $a_0 , a_{\infty}\in(p-N, p]$, $\gamma_{\infty} \leq p- a_{\infty}$ and $\gamma_{0} \geq p- a_{0}$. We first apply Theorem \ref{THM2} to get that $\lim_{R\rightarrow +\infty }\mathcal{S}_{\infty }\left(
q_{2},R\right) =0$. For this we need $q_2 > \max \left\{ 1,p\beta _{\infty },q_{*},q_{**}\right\}$, where
$$
q_{*} (\alpha_{\infty}, \beta_{\infty}, \gamma_{\infty}) =2p \frac{N+d}{2N+1}, \quad q_{**} (a_{\infty},\alpha_{\infty}, \beta_{\infty}, \gamma_{\infty}) =
p \frac{2p (d+N-1)-5}{p(2N-1) -5}.
$$

\noindent Comparing the coefficients of $d$, it is easy to see that for large $d'$s it holds $q_{**}>q_{*}$, so we obtain  $\lim_{R\rightarrow +\infty }\mathcal{R}_{\infty }\left(
q_{2},R\right) =0$ for every $q_2 > p \frac{2p (d+N-1)-5}{p(2N-1) -5} $.
To get compactness, and then existence of solutions, we now need to choose $q_1$ such that $\lim_{R\rightarrow 0^+ }\mathcal{S}_{0 }\left(
q_{1},R\right) =0$. To this end, according to Theorem \ref{THM3}, we distinguish several subcases.

\medskip

\noindent\emph{Subcase I: } $p+1 < \gamma_0 <N$. We have to choose $p<q_1 < \min \{ q_{*} (\alpha_0, \beta_0, \gamma_0),  
q_{**} (a_0 , \alpha_0, \beta_0, \gamma_0) \}$, where 
$$
q_{*} (\alpha_0, \beta_0, \gamma_0)= p \frac{2N+1}{2N-2 \gamma_0 }, \quad q_{**} (a_0 , \alpha_0, \beta_0, \gamma_0)= p \frac{\frac{p}{2}  +\gamma_0 +p(N-1) -1}
{p(N-1) -\gamma_0 (p-1)-1 }.
$$
\noindent Setting $\delta =\min \{ q_{*} (\alpha_0, \beta_0, \gamma_0),  
q_{**} (a_0 , \alpha_0, \beta_0, \gamma_0) \}$, it easy to check that $p<\delta$, and of course for large $d'$s we have $\delta <p \frac{2p (d+N-1)-3}{p(2N-1) -5} $. As a conclusion, we have compactness and thus existence of solutions for $q_1 , q_2$ such that
$$
p<q_1 < \delta < p \frac{2p (d+N-1)-3}{p(2N-1) -5} <q_2 .
$$

\medskip

\noindent\emph{Subcase II: } $ \gamma_0 =N$. In this case hypothesis $ \alpha_0 > -(1-\beta_0 )N$ is satisfied, and similar computations as above show that we can apply Theorem \ref{THM:ex2} by choosing $d$ large enough and $q_1 , q_2$ such that
$$
p<q_1 < p \frac{ \frac{p}{2} + (N-1)(p+1)    }{N-p-1}< p \frac{2p (b+N-1)-5}{p(2N-1) -5} < q_2 .
$$

\medskip

\noindent\emph{Subcase III: } $ N<\gamma_0 < \frac{p(N-1)-1}{p-1}$. In this case we have $q_{*} (\alpha_0, \beta_0, \gamma_0)<0$, and we get the result by choosing 
$$
p<q_1 < p \frac{ \frac{p}{2} + \gamma_0 +p(N-1) -1    }{p(N-1) -\gamma_0 (p-1)-1}< p \frac{2p (d+N-1)-5}{p(2N-1) -5} < q_2 .
$$

\end{exa}

\section{Appendix}

This Appendix is devoted to complete the computations of the proof of Theorem \ref{THM3}.
We still distinguish the same cases considered there.\medskip

\noindent \emph{Case }$p-a_0 <\gamma _{0}<N$.\smallskip

\noindent In this case, $\left( \alpha _{0},q_{1}\right) \in \mathcal{A}
_{a_0 ,\beta _{0},\gamma _{0}}$ means 
$$\max \left\{ 1,p\beta _{0}\right\} <q_{1}<\displaystyle\min \left\{ p\frac{\alpha
_{0}-\beta _{0}\gamma _{0}+N}{N-\gamma _{0}}, p\frac{p\alpha _{0}+\left(
1-p\beta _{0}\right) \gamma _{0}+p(N-1)+a_0}{p(N-1)-(p-1)\gamma _{0}+a_0}\right\}. $$

\noindent This ensures that we can find $\xi \geq 0$ such that 
\[
\frac{1}{p}\leq \beta _{0}+\xi \leq 1\quad \text{and}\quad p\left( \beta
_{0}+\xi \right) <q_{1}<\frac{p^2 \left( \alpha _{0}+\xi \gamma _{0}\right)
+p^2 N-p \left( \beta _{0}+\xi \right) \left( (p-1) \gamma _{0}+p-a_0 \right)  }{p(N-1)-(p-1)\gamma
_{0} +a_0}. 
\]

\noindent Indeed, this amounts to find $\xi $ such that 
\[
\left\{ 
\begin{array}{l}
\max \left\{ 0,\frac{1-p\beta _{0}}{p}\right\} \leq \xi \leq 1-\beta
_{0}\medskip \\ 
\xi <\frac{q_{1}-p\beta _{0}}{p}\medskip \\ 
q_{1}\frac{p(N-1) -\gamma_0 (p-1) +a_0 }{p(\gamma_0 -p + a_0 )}-\frac{p( \alpha _{0}+N) -\beta _{0}\left( (p-1)\gamma _{0}+p-a_0\right) }{
\gamma _{0}-p +a_0}<\xi ,
\end{array}
\right. 
\]

\noindent Since $\frac{1}{p}-\beta _{0} < 1-\beta _{0}$ is obvious (recall that $p>1$) and $1-\beta
_{0}\geq 0$ holds by assumption, such a system has a solution $\xi $ if and
only if 
\[
\left\{ 
\begin{array}{l}
\frac{1-p\beta _{0}}{p}<\frac{q_{1}-p\beta _{0}}{p}\medskip \\ 
q_{1}\frac{p(N-1)-(p-1)\gamma _{0}+a_0}{p\left( \gamma _{0}-p +a_0 \right) }-\frac{p(\alpha
_{0}+N)-\beta _{0}\left( (p-1) \gamma _{0}+p -a_0 \right) }{\gamma _{0}-p +a_0}<1-\beta
_{0}\medskip \\ 
q_{1}\frac{p(N-1)-(p-1)\gamma _{0}+a_0 }{p\left( \gamma _{0}-p +a_0\right) }-\frac{p(\alpha
_{0}+N)-\beta _{0}\left( (p-1) \gamma _{0}+p-a_0 \right) }{\gamma _{0}-p+a_0 }<\frac{
q_{1}-p\beta _{0}}{p}\medskip \\ 
\frac{q_{1}-p\beta _{0}}{p}>0,
\end{array}
\right. 
\]
which is equivalent to 
\[
\left\{ 
\begin{array}{l}
1<q_{1}\medskip \\ 
q_{1}<p\frac{p\alpha _{0}+p(N-1)+\left( 1-p\beta _{0}\right) \gamma _{0}+a_0}{
p(N-1)-(p-1)\gamma _{0}+a_0}=q_{**}\medskip \\ 
q_{1}<p\frac{\alpha _{0}+N-\gamma _{0}\beta _{0}}{N-\gamma _{0}}=q_{*}
\medskip \\ 
q_{1}>p\beta _{0}.
\end{array}
\right. 
\]

\noindent \emph{Case }$\gamma _{0}=N$.\smallskip

\noindent In this case, $\left( \alpha _{0},q_{1}\right) \in \mathcal{A}
_{a_0 ,\beta _{0},\gamma _{0}}$ means 
\[
\begin{tabular}{l}
$\alpha _{0}>\alpha _{1}\,\left( =\alpha _{2}=\alpha _{3}\right)$\quad and \smallskip\\
$\max \left\{ 1,p\beta _{0}\right\} <q_{1}<\displaystyle
p\frac{p\alpha_{0}+\left( 1-p\beta _{0}\right) \gamma _{0}+p(N-1)+a_0 }{p(N-1)-(p-1)\gamma _{0}+a_0 }=p\frac{
p\alpha _{0}+(p+1)N-p\beta _{0}N-p +a_0}{N-p+a_0}   $ 
\end{tabular}
\]
and this ensures that we can find $\xi \geq 0$ such that 
\[
\frac{1}{p}\leq \beta _{0}+\xi \leq 1\quad \text{and}\quad p\left( \beta
_{0}+\xi \right) <q_{1}<\frac{p^2 \left( \alpha _{0}+\xi \gamma _{0}\right)
+p^2 N-p\left( (p-1) \gamma _{0}+p -a_0 \right) \left( \beta _{0}+\xi \right) }{p(N-1)-(p-1) \gamma
_{0}+a_0 }. 
\]

\noindent This amounts to find $\xi $ such that 
\[
\left\{ 
\begin{array}{l}
\max \left\{ 0,\frac{1-p\beta _{0}}{p}\right\} \leq \xi \leq 1-\beta
_{0}\medskip \\ 
\frac{q_{1}}{p}-\frac{p\alpha _{0}+pN-\left((p-1) N+p-a_0 \right) \beta _{0}}{N-p+a_0 }<\xi
<\frac{q_{1}-p\beta _{0}}{p},
\end{array}
\right. 
\]
which has a solution $\xi $ if and only if 
\[
\left\{ 
\begin{array}{l}
\frac{1-p\beta _{0}}{p}<\frac{q_{1}-p\beta _{0}}{p}\medskip \\ 
\frac{q_{1}}{p}-\frac{p\alpha _{0}+pN-\left( N(p-1)+p-a_0 \right) \beta _{0}}{N-p+a_0}
<1-\beta _{0}\medskip \\ 
\frac{q_{1}}{p}-\frac{p\alpha _{0}+pN-\left( N(p-1)+p-a_0 \right) \beta _{0}}{N-p +a_0}<
\frac{q_{1}-p\beta _{0}}{p}\medskip \\ 
0<\frac{q_{1}-p\beta _{0}}{p}.
\end{array}
\right. 
\]
These conditions are equivalent to 
\[
\left\{ 
\begin{array}{l}
1<q_{1}\medskip \\ 
q_{1}<p\frac{p\alpha _0+(p+1)N-p(\beta_0 N +1) +a_0}{N-p +a_0} =q_{**}  \medskip \\ 
0<\frac{p\alpha _{0}+pN-\left( N(p-1)+p-a_0 \right) \beta _{0}}{N-p+a_0}-\beta _{0}=p\frac{
\alpha _{0}+N\left( 1-\beta _{0}\right) }{N-p+a_0}=p\frac{\alpha _{0}-\alpha _{1}
}{N-p+a_0}\medskip \\ 
p\beta _{0}<q_{1}.
\end{array}
\right. 
\]

\noindent \emph{Case }$N<\gamma _{0}<\frac{p(N-1)+a_0}{p-1}$.\smallskip 

\noindent In this case, $\left( \alpha _{0},q_{1}\right) \in \mathcal{A}
_{a_0 ,\beta _{0},\gamma _{0}}$ means 
\[
\begin{tabular}{l}
$\displaystyle\max \left\{ 1,p\beta _{0},p\frac{\alpha _{0}-\beta _{0}\gamma _{0}+N}{N-\gamma _{0}}\right\} <q_{1}<p%
\frac{p\alpha _{0}+\left( 1-p\beta _{0}\right) \gamma _{0}+p(N-1)+a_0}{p(N-1)-(p-1)\gamma_{0}+a_0}$
\end{tabular}
\]
and these conditions ensure that we can find $\xi \geq 0$ such that 
\[
\frac{1}{p}\leq \beta _{0}+\xi \leq 1\quad \text{and}\quad p\left( \beta
_{0}+\xi \right) <q_{1}<\frac{p^2 \left( \alpha _{0}+\xi \gamma _{0}\right)
+p^2 N-p\left( (p-1) \gamma _{0}+p -a_0\right) \left( \beta _{0}+\xi \right) }{(N-1)p-(p-1)\gamma
_{0}+a_0}. 
\]
Indeed, this is equivalent to find $\xi $ such that 
\[
\left\{ 
\begin{array}{l}
\max \left\{ 0,\frac{1-p\beta _{0}}{p}\right\} \leq \xi \leq 1-\beta
_{0}\medskip \\ 
\xi <\frac{q_{1}-p\beta _{0}}{p}\medskip \\ 
q_{1} \frac{p(N-1) -\gamma_0 (p-1) +a_0}{p(\gamma_0 -p +a_0 )}    -\frac{p \alpha _{0}+p N-\beta _{0}\left( (p-1)\gamma _{0}+p-a_0\right) }
{\gamma _{0}-p+a_0}<\xi ,
\end{array}
\right. 
\]
Such a system has a solution $\xi $ if and only if 
\[
\left\{ 
\begin{array}{l}
0<\frac{q_{1}-p\beta _{0}}{p}\medskip \\ 
\frac{1-p\beta _{0}}{p}<\frac{q_{1}-p\beta _{0}}{p}\medskip \\ 
q_{1}\frac{p(N-1)-(p-1)\gamma _{0}+a_0}{p\left( \gamma _{0}-p +a_0\right) }-\frac{p\alpha
_{0}+pN-\left( (p-1)\gamma _{0}+p -a_0\right) \beta _{0}}{\gamma _{0}-p +a_0}  <1-\beta
_{0}\medskip \\ 
q_{1}\frac{p(N-1)-(p-1)\gamma _{0} +a_0}{p\left( \gamma _{0}-p +a_0\right) }-\frac{p\alpha
_{0}+pN-\left( (p-1)\gamma _{0}+p -a_0 \right) \beta _{0}}{\gamma _{0}-p +a_0} <\frac{q_{1}-p\beta _{0}}{p},
\end{array}
\right. 
\]
which is equivalent to 
\[
\left\{ 
\begin{array}{l}
q_{1}>p\beta _{0}\medskip \\ 
q_{1}>1\medskip \\ 
q_{1}<p\frac{p\alpha _{0}+p(N-1)+ \gamma_0 (1-p\beta _{0}) +a_0}{
p(N-1)-(p-1)\gamma _{0}+a_0}=q_{**}\medskip \\ 
q_{1}> p\frac{\alpha _{0}+N-\beta _{0}\gamma _{0}}{N-\gamma _{0}}= q_{*}
\end{array}
\right. 
\]

\noindent \emph{Case }$\gamma _{0}=\frac{p(N-1)+a_0}{p-1}$\emph{.\smallskip }

\noindent In this case, $\left( \alpha _{0},q_{1}\right) \in \mathcal{A}
_{a_0 , \beta _{0},\gamma _{0}}$ means 
\[
\alpha _{0}>\alpha _{1}\quad \text{and}\quad q_{1}>\max \left\{ 1,p\beta
_{0},-p\frac{(\alpha _{0}+N)(p-1)-\beta_{0}(p(N-1)+a_0 )}{N-p +a_0}\right\} . 
\]
This ensures that we can find $\xi \geq 0$ such that 
$$
\frac{1}{p}\leq \beta _{0}+\xi \leq 1,\quad q_{1}>p\left( \beta _{0}+\xi
\right)$$
\noindent and
$$p \left( \alpha_0 + \xi \gamma_0 \right) +pN  - \left( \left( p-1 \right) \gamma_0 +p -a_0 \right) \left( \beta_0 + \xi \right)>0, 
$$
\noindent that is 
\[
\frac{1}{p}-\beta _{0}\leq \xi \leq 1-\beta _{0},\quad q_{1}>p\beta
_{0}+p\xi \quad \text{and}\quad \frac{p}{p-1} \left(  N-p+a_0  \right) \xi +p \left(  \alpha_0 +N \left( 1-\beta_0 \right)  \right)>0. 
\]
Indeed, this amounts to find $\xi $ such that 
\[
\left\{ 
\begin{array}{l}
\max \left\{ 0,\frac{1-p\beta _{0}}{p}\right\} \leq \xi \leq 1-\beta
_{0}\medskip \\ 
-\frac{(p-1) \left(\alpha _{0}+N\left( 1-\beta _{0}\right) \right)}{N-p+a_0 }<\xi <\frac{
q_{1}-p\beta _{0}}{p}
\end{array}
\right. 
\]
and such a system has a solution $\xi $ if and only if 
\[
\left\{ 
\begin{array}{l}
\frac{1-p\beta _{0}}{p}<\frac{q_{1}-p\beta _{0}}{p}\medskip \\ 
-\frac{(p-1) \left( \alpha _{0}+N\left( 1-\beta _{0}\right) \right)  }{N-p+a_0} <1-\beta _{0}\medskip
\\ 
-\frac{(p-1) \left( \alpha _{0}+N\left( 1-\beta _{0}\right) \right)  }{N-p+a_0}  <\frac{q_{1}-p\beta _{0}
}{p}\medskip \\ 
0<\frac{q_{1}-p\beta _{0}}{p},
\end{array}
\right. 
\]
which means 
\[
\left\{ 
\begin{array}{l}
1<q_{1}\medskip \\ 
\alpha _{0}>-\left( 1-\beta _{0}\right) \frac{p(N-1) +a_0}{p-1}  =-\left( 1-\beta _{0}\right) \gamma_0 = \alpha
_{1}\medskip \\ 
q_{1}>p\left( \beta _{0}-\frac{(p-1) \left( \alpha _{0}+N\left( 1-\beta _{0}\right) \right) }{N-p+a_0} \right)
=-p\frac{ (p-1)\alpha_0 -\beta_0 p (N-1) -\beta_0 a_0  +N(p-1)}{N-p+a_0}=q_{*}\medskip \\ 
p\beta _{0}<q_{1}.
\end{array}
\right. 
\]

\noindent \emph{Case }$\gamma _{0}>\frac{p(N-1)+a_0}{p-1} $.\smallskip

\noindent In this case, $\left( \alpha _{0},q_{1}\right) \in \mathcal{A}
_{a_0 , \beta _{0},\gamma _{0}}$ means 
\[
q_{1}>\max \left\{ 1,p\beta _{0},p\frac{\alpha _{0}-\beta _{0}\gamma _{0}+N}{
N-\gamma _{0}},p\frac{p\alpha _{0}+\left( 1-p\beta _{0}\right) \gamma
_{0}+p(N-1)+a_0}{p(N-1)-(p-1)\gamma _{0}+a_0}\right\} 
\]
and this condition ensures that we can find $\xi \geq 0$ such that 
\[
\frac{1}{p}\leq \beta _{0}+\xi \leq 1\quad \text{and}\quad q_{1}>p\max
\left\{ \beta _{0}+\xi ,\frac{p\left( \alpha _{0}+\xi \gamma _{0}\right)
+pN-\left( (p-1)\gamma _{0}+p-a_0 \right) \left( \beta _{0}+\xi \right) }{p(N-1)-(p-1)\gamma
_{0}+a_0 }\right\} , 
\]
which amounts to find $\xi $ such that 
\begin{equation}
\left\{ 
\begin{array}{l}
\max \left\{ 0,\frac{1-p\beta _{0}}{p}\right\} \leq \xi \leq 1-\beta
_{0}\medskip \\ 
\frac{q_{1}}{p}>\max \left\{ \beta _{0}+\xi ,\frac{\gamma _{0}-p+a_0}{
p(N-1)-(p-1)\gamma _{0}+a_0 }\xi +\frac{p\alpha _{0}+pN-\beta _{0}\left( (p-1)\gamma _{0}+p-a_0 \right)
}{p(N-1)-(p-1)\gamma _{0}+a_0}\right\} .
\end{array}
\right.  \label{ineq}
\end{equation}
In order to check this, we take into account that $\gamma _{0}>\frac{p(N-1)+a_0}{p-1}$ and $a_0 >p-N$ imply 
$\gamma _{0}>N$, and observe that 
\[
\beta _{0}+\xi =\frac{\gamma _{0}-p+a_0}{p(N-1)-(p-1)\gamma _{0}+a_0}\xi +\frac{p\alpha
_{0}+pN-\beta _{0}\left( (p-1)\gamma _{0}+p-a_0\right) }{p(N-1)-(p-1)\gamma _{0}+a_0}
\Longleftrightarrow \xi =\frac{\alpha _{0}+\left( 1-\beta _{0}\right) N}{
N-\gamma _{0}}. 
\]
Accordingly, we distinguish three subcases: \medskip 

\noindent (I) $\frac{\alpha _{0}+\left( 1-\beta _{0}\right) N}{N-\gamma _{0}}
\geq 1-\beta _{0}$, i.e., $\alpha _{0}\leq -\gamma _{0}\left( 1-\beta
_{0}\right) =\alpha _{1}$;$\medskip $

\noindent (II) $\frac{\alpha _{0}+\left( 1-\beta _{0}\right) N}{N-\gamma _{0}%
}\leq \max \left\{ 0,\frac{1-p\beta _{0}}{p}\right\} $, i.e., 
\[
\alpha _{0}+\left( 1-\beta _{0}\right) N\geq \left( N-\gamma _{0}\right)
\max \left\{ 0,\frac{1-p\beta _{0}}{p}\right\} =\min \left\{ 0,\left(
N-\gamma _{0}\right) \frac{1-p\beta _{0}}{p}\right\} , 
\]
i.e., 
\[
\alpha _{0}\geq \min \left\{ 0,\left( N-\gamma _{0}\right) \frac{1-p\beta
_{0}}{p}\right\} -\left( 1-\beta _{0}\right) N=\min \left\{ \alpha
_{2},\alpha _{3}\right\} ; 
\]

\noindent (III) $\max \left\{ 0,\frac{1-p\beta _{0}}{p}\right\} <\frac{%
\alpha _{0}+\left( 1-\beta _{0}\right) N}{N-\gamma _{0}}<1-\beta _{0}$,
i.e., $\alpha _{1}<\alpha _{0}<\min \left\{ \alpha _{2},\alpha _{3}\right\}
.\medskip $

\noindent \emph{Subcase (I).\smallskip }

\noindent Since $\xi \leq 1-\beta _{0}$ implies 
\begin{eqnarray*}
&& 
\max \left\{ \beta _{0}+\xi ,\frac{\gamma _{0}-p+a_0 }{p(N-1)-(p-1)\gamma _{0}+a_0}\xi 
+ \frac{p\alpha _{0}+pN-\beta _{0}\left( (p-1)\gamma _{0}+p -a_0 \right) }{p(N-1)-(p-1)\gamma_{0}+a_0} \right\} \\
& = &
\frac{\gamma _{0}-p+a_0 }{p(N-1)-(p-1)\gamma _{0}+a_0}\xi 
+ \frac{p\alpha _{0}+pN-\beta _{0}\left( (p-1)\gamma _{0}+p -a_0 \right) }{p(N-1)-(p-1)\gamma_{0}+a_0} ,
\end{eqnarray*}
the inequalities (\ref{ineq}) become 
\[
\left\{ 
\begin{array}{l}
\max \left\{ 0,\frac{1-p\beta _{0}}{p}\right\} \leq \xi \leq 1-\beta
_{0}\medskip \\ 
\frac{q_{1}}{p}>\frac{\gamma _{0}-p+a_0 }{p(N-1)-(p-1)\gamma _{0}+a_0}\xi 
+ \frac{p\alpha _{0}+pN-\beta _{0}\left( (p-1)\gamma _{0}+p -a_0 \right) }{p(N-1)-(p-1)\gamma_{0}+a_0} ,
\end{array}
\right. 
\]
i.e., 
\[
\left\{ 
\begin{array}{l}
\max \left\{ 0,\frac{1-p\beta _{0}}{p}\right\} \leq \xi \leq 1-\beta
_{0}\medskip \\ 
q_{1}\frac{p(N-1)-(p-1)\gamma _{0}+a_0}{p\left( \gamma _{0}-p+a_0\right) }-\frac{p^2 \alpha
_{0}+p^2 N-p \beta _{0}\left( (p-1)\gamma _{0}+p-a_0\right) }{p\left( \gamma
_{0}-p+a_0\right) } <\xi ,
\end{array}
\right. 
\]
which, since $\max \left\{ 0,\frac{1-p\beta _{0}}{p}\right\} \leq 1-\beta
_{0}$ is clearly true, has a solution $\xi $ if and only if 
\[
q_{1}\frac{p(N-1)-(p-1)\gamma _{0}+a_0}{p\left( \gamma _{0}-p+a_0 \right) }-\frac{p^2 \alpha
_{0}+p^2 N-p\beta _{0}\left( (p-1)\gamma _{0}+p-a_0\right) }{p\left( \gamma
_{0}-p+a_0\right) }  <1-\beta _{0}, 
\]
i.e., 
\begin{eqnarray*}
q_{1}&>&
p\frac{p\alpha _{0}+pN-p+\gamma _{0}\left( 1-p\beta _{0}\right) +a_0}{p(N-1)-(p-1)\gamma _{0}+a_0}=q_{**}. 
\end{eqnarray*}

\noindent \emph{Subcase (II).\smallskip }

\noindent Since $\xi \geq \max \left\{ 0,\frac{1-p\beta _{0}}{p}\right\} $
implies $\max \left\{ \beta _{0}+\xi ,\frac{\gamma _{0}-p+a_0}{p(N-1)-(p-1)\gamma _{0}+a_0}
\xi +\frac{p\alpha _{0}+pN-\beta _{0}\left( (p-1)\gamma _{0}+p-a_0\right) }{
p(N-1)-(p-1)\gamma _{0}+a_0}\right\} =\beta _{0}+\xi $, the inequalities (\ref{ineq})
become 
\[
\left\{ 
\begin{array}{l}
\max \left\{ 0,\frac{1-p\beta _{0}}{p}\right\} \leq \xi \leq 1-\beta
_{0}\medskip \\ 
\xi <\frac{q_{1}}{p}-\beta _{0},
\end{array}
\right. 
\]
which has a solution $\xi $ if and only if $\max \left\{ 0,\frac{1-p\beta
_{0}}{p}\right\} \leq \frac{q_{1}}{p}-\beta _{0}$, i.e., $q_{1}>\max \left\{
1,p\beta _{0}\right\} $.$\medskip $

\noindent \emph{Subcase (III).\smallskip }

\noindent We take $\xi =\frac{\alpha _{0}+\left( 1-\beta _{0}\right) N}{
N-\gamma _{0}}$ and thus the inequalities (\ref{ineq}) are equivalent just
to 
\[
\frac{q_{1}}{p}>\beta _{0}+\frac{\alpha _{0}+\left( 1-\beta _{0}\right) N}{
N-\gamma _{0}}=\frac{\alpha _{0}+N-\gamma _{0}\beta _{0}}{N-\gamma _{0}}, 
\]
that is $ q_1 > q_{*}$.



\end{document}